\documentclass{amsart}
\usepackage{amssymb}
\usepackage{amscd}
\usepackage{amsmath}
\usepackage{amsfonts}
\usepackage{a4}
\usepackage{amsthm}

\input xy
\xyoption{all}

\newtheorem{thm}{Theorem}[section]
\newtheorem{defi}[thm]{Definition}
\newtheorem{prop}[thm]{Proposition}
\newtheorem{lem}[thm]{Lemma}
\newtheorem{cor}[thm]{Corollary}

\newenvironment{rem}{\noindent \bf Remark : \rm \newline \noindent}{}

\newcommand{\ra}{\rightarrow}
\newcommand{\thra}{\twoheadrightarrow}

\newcommand{\mc}{\mathcal }
\newcommand{\udl}{\underline }
\newcommand{\ab}{{\rm \bf Ab} }
\renewcommand{\mod}{{\rm \bf mod} \,}
\newcommand{\Mod}{{\rm \bf Mod} \,}
\newcommand{\rad}{{\rm rad} \,}
\newcommand{\spec}{{\rm Spec} \,}
\newcommand{\Qcoh}{{\rm \bf Qcoh} \,}
\newcommand{\qcoh}{{\rm  Qcoh} \,}
\newcommand{\Coh}{{\rm \bf Coh} \,}
\newcommand{\coh}{{\rm  Coh} \,}

\begin{document}

\bibliographystyle{plain}

\author{Niels Borne}
\address{Niels Borne,
Universit\`a di Bologna \\
Dipartimento di Matematica \\
Piazza di Porta San Donato, 5 \\
I-40126 Bologna \\
Italy}
\email{borne@dm.unibo.it}
\title[Cohomology of $G$-sheaves in positive characteristic]
{Cohomology of $G$-sheaves in positive characteristic}
\maketitle
\tableofcontents


\section{Introduction}

\subsection{Galois modules in positive characteristic}

This work is an attempt to a better understanding of Galois modules in
positive characteristic. We define them here as spaces of global sections
of coherent $G$-sheaves on a proper scheme $X$ over an algebraically
closed field $k$ of positive characteristic,
endowed with the action of a finite group $G$.
So Galois modules are representations of $G$, and they are known to
be related to the ramification of the action (i.e. to the fixed points).

When the characteristic $p$ of $k$ does not divide the order of
$G$, i.e. in the reductive case, Galois modules are quite well
understood, and the situation is similar to the one in
characteristic zero. But when $p$ does divide the order of $G$,
the situation is much more mysterious.

\subsection{The role of equivariant $K$-theory}
Indeed, the general approach to the description of Galois modules
is the use of equivariant $K$-theory, which provides a
Euler-Poincar\'e characteristic $\chi(G,\cdot) : K_0(G, X)\ra
R_k(G)$ with values in the Brauer characters group of $G$. One can
use an equivariant Lefschetz formula to compute it explicitly.
This is satisfactory in the reductive case, since then the Brauer
character of a representation characterizes its isomorphism class.
But this last fact is false as soon as $p$ divides the order of
$G$, because then the Galois modules are modular representations,
so their Brauer characters only describe their Jordan-H\"older
series, and not their decomposition in indecomposables (or
equivalently their isomorphism classes). This is particulary
dramatic when $G$ is a $p$-group, because then
$R_k(G)\simeq\mathbb Z$ contains no equivariant information at
all.

The usual way to improve the situation is to make assumptions on
the ramification of the action. For an example, on a smooth
projective curve, it one assumes that the action is tamely
ramified, one can show that a coherent $G$-sheaf of large degree
has a space of global sections  which is a projective
$k[G]$-module, so that its Brauer character is enough to describe
its isomorphism class. This is part of the Noether-like criterion
given by S.Nakajima in \cite{Naka}, and this point of view was
developed later on (see \cite{ChinEr}).

For a wild action, equivariant $K$-theory is much less effective.
In particular, it is unable to explain the
modular Riemann-Roch formula given by S.Nakajima in \cite{Naka}, which
describes the structure of the space of global sections $H^0(X,
\mc L)$ of an invertible sheaf of large degree on a smooth projective
curve with an action of $G=\mathbb Z /p$.

A first look at the properties a refined version of equivariant
$K$-theory should check suggests that important changes are
needed. Indeed, the Lefschetz formula shows that if the action of
$G$ is free, the equivariant Euler-Poincar\'e characteristic is a
multiple of $[k[G]]$ in $R_k(G)$. However, another work of
S.Nakajima (see \cite{NakaInv}) shows that such a ``symmetry
principle'' fails to hold in a refined sense if we use the
standard Zariski cohomology of sheaves to define the refined
Euler-Poincar\'e characteristic.

\subsection{Modular $K$-theory : main properties}

In this work, we present a refinement of the equivariant $K$-theory of
a noetherian $G$-scheme $X$ over $k$. To do so, we introduce,
for each full subcategory $\mc A$ of the category $k[G]\,\mod$ of
$k[G]$-modules of finite type, the notion of $\mc A$-sheaf on $X$.
This is, roughly, a sheaf of
modules over the Auslander algebra $\mc A_X$, itself defined as a
certain sheaf of algebras with several objects over the quotient scheme
$Y=X/G$, obtained by mimicking the functorial
definition of the Auslander algebra of $k[G]$ (see \cite{Aus}).

We imbed
the category $\qcoh(G,X)$ of quasicoherent $G$-sheaves on $X$
as a reflective subcategory of the category
$\qcoh(\mc A,X)$ of quasicoherent $\mc A$-sheaves on $X$
(see \S \ref{Gsh>Ash}), so that each
$G$-sheaf $\mc F$ can be seen as a $\mc A$-sheaf $\udl{\mc F}$.
The category of coherent $\mc A$-sheaves on $X$ is abelian, and its
$K$-theory, in the sense of Quillen, is denoted by $K_i(\mc A,X)$.
This construction is functorial in both variables.

Of course, the main test of validity for this new definition is the
case $X= {\rm spec}\, k$. Modular $K$-theory is
satisfying for groups with cyclic $p$-Sylows, where $p$ is the
characteristic of $k$, in the sense that for each $k[G]$-module of
finite type $V$, its class $[\udl {V}]$ in $K_0(\mc A, {\rm spec}\, k)$
characterizes its isomorphism class (see Theorem \ref{Main}).
For an arbitrary finite group $G$,
one grasps certainly more information than the Brauer character,
but until now it is not clear to the author if the definition enables
to get back the isomorphism class.

Using the functoriality in $\mc A$, we can compare both $K$-theories.
Indeed if $\mc A$ and $\mc A'$, seen as rings with several objects,
are Morita equivalent,
then their modules and their $K$-theory are the same (see Proposition
\ref{CohMor}). In particular,
since the category $\mc P$ of projective $k[G]$-modules of finite type
is Morita equivalent to the category with only the free object $k[G]$,
the groups $K_i(\mc P, X)$ coincide with the equivariant
$K$-theory. So in the reductive case, since $\mc P=\mc A$,
we get nothing else than equivariant $K$-theory.
Hence something new happens only when we consider a $\mc A$ checking $\mc
P \subsetneqq \mc A$, and in this case we have
a surjective homomorphism $K_0(\mc A, X) \thra K_0(G, X)$.

Since $K_i(\mc A, X)$ has the usual functorial properties in $X$ (or
rather in the quotient $Y=X/G$), it makes sense to ask whether there is a
localization long exact sequence for modular $K$-theory. We answer
positively (see Theorem \ref{Loc}), but only on a certain
surjectivity assumption, which is always fulfilled for equivariant
$K$-theory, but not for modular $K$-theory.

The computation of modular Euler-Poincar\'e characteristics
reduces to the computation of standard Euler-Poincar\'e
characteristics on the quotient scheme,thanks to the formula~:

$$\chi(\mc A, \mc F)=\sum_{I\in S}\chi(\mc F(I))[S_I]$$
(see Lemma \ref{red}). Of course, the modular Euler-Poincar\'e
characteristic here is computed with the cohomology of $\udl{\mc F}$,
which is very different from the cohomology of $\mc F$. We give
comparison results in section \ref{Cohomology}.

Among the applications, we show that the symmetry principle holds in
modular $K$-theory. Indeed, if the action of $G$ on $X$ is free, then
for any $G$-sheaf $\mc F$ we have an equality in $K_0(\mc
A, {\rm spec}\, k)$ :

$$\chi(\mc A, \udl {\mc F})=\chi(\pi_*^G \mc F) \left[ \udl{k[G]}
\right] $$
(see Proposition \ref{sym}).

The other applications all concern the $1$-dimensional case, i.e.
the case when $X$ is a projective curve over $k$.

Thanks to the localization sequence, we give a description
of the additive structure of $K_0(X, \mc A)$ (see Theorem \ref{Mad})
when the group acts with normal stabilizers.
This is done by introducing a group of class of cycles with
coefficients in the modular representations $A_0(\mc A,X)$, and enables
to define a first Chern class such that the usual Riemann-Roch formula
holds. We do not compute these modular Chern classes in general
(although it is likely possible) but only those for invertible sheaves
pulled back from the quotient, which have a very simple expression. As
a reward, we get :

\begin{thm}

Let $X$ be a projective curve over an algebraically closed field $k$
of positive characteristic $p$, endowed with the faithful action of
finite group $G$. Suppose that $G$ has cyclic $p$-Sylows  and
acts with normal stabilizers.
Denote by $\pi : X \ra Y=X/G$ the quotient map.
There exists an integer $N$ such that for all invertible sheaves $\mc
M$, $\mc N$ on $Y$ such that $\deg \mc M \geq \deg  \mc N \geq N$ the
following isomorphism of $k[G]$-modules holds :

$$ H^0(X, \pi^* \mc M) \simeq H^0(X, \pi^*\mc N) \oplus k[G]^{\oplus
  \deg \mc M - \deg \mc N}$$

\end{thm}

\noindent (see Theorem \ref{pullbackinv}).

To end with, we extend the result of S.Nakajima in \cite{Naka} to
the case of an action of a cyclic group, giving a recursive
algorithm to compute explicitly the structure of modular
representation of the space $H^0(X, \mc L)$ of global sections of
an invertible sheaf of large degree on the curve (see section
\ref{Cyclic}). In this case the modular Euler-Poincar\'e
characteristics takes the explicit form :

$$ \chi (\mc A, \udl {\mc L})
=\sum_{\psi \in \widehat G} \sum_{j=1}^{p^v} \sum_{i=1}^{j}
\chi({\rm gr}_0 \udl {\mc L}(\psi \otimes V_{i}))[S_{\psi \otimes
V_{j}}]$$ (see Lemma \ref{explicit}). The ${\rm gr}_0 \udl {\mc
L}(\psi \otimes V_{j})$ are invertible sheaves on $Y$, and can be
represented by divisors thanks to the key Proposition \ref{recur}
which describes their behaviour under extension :

$$ {\rm gr}_0^G  \udl {\mc L}(V_j) \simeq
{\rm gr}_0^P \udl{{\rm gr}_0^H\udl{\mc L}(V_l)}(V_{j'})$$
and thus allows to reduce to the case of a cyclic group of order $p$.

To illustrate the use of our algorithm, we give an explicit
expression of the structure of $H^0(X, \mc L)$ when $X$ is a
cyclic $p$-group : see Theorem \ref{formula}. The expression of
the coefficient of a given indecomposable in this decomposition
involves the use of all the ramification jumps of the
corresponding cover.

From this rather brutal computation we can deduce a Noether-like
criterion, showing in a more qualitative way how ramification and
Galois modules are linked (see Theorem \ref{Noether}).
Remember that a $k[G]$-module is said \emph{relatively
$H$-projective} if it is a direct summand of a module induced from
$H$.

\begin{thm}
Let $\pi : X\ra Y$ be a (generically) cyclic Galois $p$-cover of
projective curves over $k$ of group $G$, ${\rm ram}\,\pi$ the
largest ramification subgroup of $\pi$, and $H$ a subgroup of $G$.
Then the following assertions are equivalent :

(i) ${\rm ram}\,\pi \subset H$

(ii) $\forall \mc L \in {\rm Pic}_G X \;\;\; \deg \mc L > 2g_X-2
\Longrightarrow H^0(X, \mc L)$ is relatively $H$-projective.

(iii) $\exists \mc M \in {\rm Pic}\, Y$ so that $\# G \deg \mc M >
2g_X-2$ and $H^0(X, \pi^*\mc M)$ is relatively $H$-projective.

\end{thm}

Even if it is likely that a direct cohomological proof exists, the
author was unable to find one.

\subsection{Modular $K$-theory : construction}

We give here some indications about the organization of this article.

Section \ref{Repr} is devoted to the analysis of the
$0$-dimensional case. It relies in an essential way on the work of
Auslander, who first realized the interest of rings with several
objects for modular representation theory. The main idea is that
for $\mc A= k[G] \,\mod$, the isomorphism classes of
indecomposable objects of $\mc A$ are in one-to-one
correspondence, via the Yoneda embedding, with isomorphism classes
of simple objects of $\mod \mc A$, and so are well detected by
$K$-theory there. Only the reinterpretation in terms of
Grothendieck groups we give seems not to have been used before,
even if it is not surprising in itself. The strongest result,
obtained for groups with cyclic $p$-Sylows, is Theorem \ref{Main}.

In section \ref{Kan}, we sum up briefly the tools of enriched category
theory needed in the sequel. We give a proof only of the
facts that we could not find in the literature.

Section \ref{RingedSchemes} is a preliminary to the next section on
$\mc A$-sheaves. We define a ringed scheme as a scheme $Y$, endowed
with a category enriched in the closed category $\Qcoh Y$,
which is a way to express the notion of
``scheme with a sheaf of algebras with several objects''. The
Auslander algebra $\mc A_X$ of a $G$-scheme $X$ will be an
example. However, it seems convenient to deal with some problems at this
level of generality, especially the problem of showing the existence
of an adjunction between pull-back and push-forward for sheaves of
modules over ringed schemes (see Proposition \ref{AdjunctionRingedSchemes}).

In a next step (section \ref{AuslanderAlgebra}), we define the
Auslander algebra $\mc A_X$ of a $G$-scheme $X$, which is a category
enriched on  $\Qcoh Y$, where $Y=X/G$, and introduce some variants.
Then $\mc A$-sheaves on $X$ are defined as sheaves for the ringed
scheme $(Y, \mc A_X)$, and are studied in the rest of the section.
Here, the flexibility we have in choosing the subcategory $\mc A$ of
$k[G]\, \mod$ proves useful, since additional structures on $\mc A$
give rise to additional structures on $\Qcoh(A,X)$. For instance if
$\mc A$ is monoidal, then $\Qcoh(A,X)$ is closed (see \S \ref{Amonoidal}),
as B.Day noticed in a more general context (see \cite{Day}). Another
example can be found in \ref{CanFiltr}, where it is shown that if $\mc
A$ is stable by radical and duality, and both operations commute,
then every monomorphism preserving $\mc A$-sheaf (and in particular
those coming from $G$-sheaves) has a canonical filtration,
a helpful fact to compute modular Euler-Poincar\'e characteristics.

\subsection{Acknowledgments}

This paper was written in the University of Bologna, where the author
holds a postdoctoral position. My first thanks go to A.Vistoli~:
without his support, my initial project would certainly have remained as
such for some years, if not for ever. I am also grateful to G.Vezzosi
for many enriching conversations. Finally, I thank B.K\"ock whose
comments have helped me a lot to improve on the first version.


\section{Modular representation theory following Auslander}
\label{Repr}

\subsection{Modules over a ring with several objects}

The idea that a small additive category behaves like a (non
necessarily commutative) ring is due to Mitchell (see \cite{Mitch}).
We present quickly the basic notions of the theory, following the
exposition of Auslander in \cite{Aus} (see also \cite{Street}).
We will most of the time omit the word ``small'',
although it is logically necessary in the definition of a ring with
several objects, to avoid the discussion on universes needed to make
the definitions coherent.

\subsubsection{Definition}

As usual $\ab$ denotes the category of abelian groups. All categories,
functors considered in part \ref{Repr}
are additive\footnote{We don't use the terminology
``preadditive category''. So an additive category does not need to have
finite coproducts.}, i.e. enriched over
$\ab$

For a category $\mc A$, we will denote by $\Mod \mc A$
the category $[\mc A^{op},\ab]$ of contravariant
functors from $\mc A$ to $\ab$, and natural transformations between
them. A {\it (right) $\mc A$-module} is by definition an object in $\Mod
\mc A$. $\mc A$-modules obviously form an abelian category.

The usual Yoneda embedding $\mc A \ra \Mod \mc
A$ sends the object $V$ to the contravariant representable functor
$\udl V=\mc A(\cdot,V)$. A $\mc A$-module is said {\it of finite type} if
it is a quotient of a finite direct sum of representable functors.
We will call $\mod \mc A$
the full subcategory of $\Mod \mc A$ consisting of objects of
finite type.

The {\it projective completion} of $\mc A$, denoted by $Q\mc A$, is the full
subcategory of $\Mod \mc A$ whose objects are the projective $\mc
A$-modules of finite type (equivalently, direct summands of finite
direct sums of representable functors). Since representable functors are
projective, the Yoneda embedding factorizes trough $Q\mc A$. The
resulting embedding $\mc A \ra Q\mc A$ is an equivalence if and only
if $\mc A$ is a category with finite coproducts where idempotents
split (i.e. for all $e:V \ra V$  idempotent in $\mc A$, $e$ has a
kernel in $\mc A$).

\subsubsection{Morita equivalence}
\label{Morita}
By definition, two categories $\mc A$, $\mc A'$ are {\it Morita equivalent}
if the corresponding categories of modules $\Mod \mc A$, $\Mod
\mc A'$ are equivalent. This is well known to be the case if and
only if the projective completions $Q\mc A$ and $Q\mc A'$ are
equivalent (see \cite {Aus}, Proposition 2.6).
In particular, $\mc A$ and $Q\mc A$ are Morita
equivalent.

\subsubsection{Change of ring}

\begin{prop}
\label{RightAdj}
Let $\mc A$ be a category and $\mc A'$ be a full subcategory
of $\mc A$.
The restriction functor $R:{\Mod \mc A}
\ra {\Mod \mc A'} $ admits as a right adjoint the functor $K$
defined on objects by :

\begin{equation*}
 \xymatrix @R=0pt {
  K:{\Mod \mc A'} \ar[r]  &  {\Mod \mc A} & \\
          F        \ar[r] &  (V \ar[r] & {\Mod \mc A'}
(\mc A(\cdot,V)_{|\mc A' },F))}
\end{equation*}
Moreover the counit of this adjunction is an isomorphism :
$RK \simeq 1$ (equivalently, $K$ is fully faithful).

In particular, if $\mc A$ is projectively complete,
and every object of $\mc A$ is a direct summand of a
finite direct sum of objects of $\mc A'$, this adjunction is an
equivalence ${\Mod \mc A}\simeq {\Mod \mc A'} $.
\end{prop}

\begin{proof}
See \cite{Aus}, Proposition 3.4 and Proposition 2.3.
\end{proof}

This Proposition is an example of {\it enriched Kan extension}, notion that
we will use later in a larger context : see \S \ref{Kan}. We will be
particulary interested in the following situation.

\begin{defi}
When $\mc A$ is Morita equivalent to (the full subcategory generated
by) a finite set of its objects, we will say that $\mc A$ \emph{admits
  a finite set of additive generators}.
\end{defi}

\subsection{The finite representation type case}

In this paragraph, we give an interpretation in terms of Grothendieck
groups of the classical link between
rings with several objects and modular representation theory (i.e. the study
of representations of a finite group $G$ over an
algebraically closed field $k$ of characteristic dividing the order of
$G$).

\subsubsection{The ring with several objects $\mod \mc A_{tot}$}

\label{usualeq}

More generally fix an algebraically closed field $k$,
$R$ a (non necessarily commutative) finite dimensional $k$-algebra.
Our main interest lies in the case $R=k[G]^{op}$, the opposite of
the group algebra of a finite group $G$, when the characteristic
$p$ of $k$ divides the order of the cardinal of $G$.

The idea is that to study the category $\mc A_{tot}=\mod R$, the category of
right $R$-modules of finite type, it is useful to consider it as a
ring with several objects, and thus consider the associated module
category $\Mod \mc A_{tot}$.

\begin{defi}
The algebra $R$ is said \emph{ of finite representation type} if the
number of isomorphism classes of indecomposable objects of $\mc A_{tot}$ is
finite.
\end{defi}

So it means precisely that $\mc A_{tot}$ has an additive generator.
We try to keep the letter $\mc A$ for an arbitrary subcategory of
$\mod R$ (often supposed to possess a finite set of generators),
and sometimes for an arbitrary additive category.

\subsubsection{Grothendieck groups of categories of modules}

For an arbitrary exact category $\mc A$ (see \cite{Quil}), we will
use the traditional notation $K_0(\mc A)$ to denote its
\emph{Grothendieck group}~: explicitly, it is the quotient of the
free abelian group generated by the isomorphism classes $[V]$ of
objects $V$ of $\mc A$ by the subgroup generated by the
expressions $[V]=[V']+[V'']$ associated to exact sequences $0\ra
V' \ra V \ra V'' \ra 0$ in $\mc A$. Any abelian category $\mc A$
will be endowed with its canonical exact structure. Moreover, any
additive category $\mc A$ has also a canonical exact structure
consisting of the only \emph{split} exact sequences, and we will
always use the notation $\mc A^{split}$ for it.

Since we wish to introduce finer invariants of the $R$-modules than
the usual Brauer character (i.e. for a $R$-module $V$, simply its
class $[V]$ in $K_0(\mc A_{tot})$), we have to introduce
larger groups. The immediate idea is to suppress relations in the
definition of $K_0(\mc A_{tot})$ and thus consider the group
$K_0(\mc A_{tot}^{split})$. We have an obvious epimorphism
$K_0(\mc A_{tot}^{split}) \thra K_0(\mc A_{tot})$,
and moreover the Krull-Schmidt
theorem tells us that the class $[V]$ of an object $V$ of
$\mc A_{tot}$ in $K_0(\mc A_{tot}^{split})$
determines its isomorphism class (see
Lemma \ref{KS}).

However, if the group $K_0(\mc A_{tot}^{split})$ we have just built is certainly the
right one, the way we have built it is wrong. Indeed, if we try to
export this construction in higher dimension, by considering instead
of $\mc A_{tot}$ the category $\coh(G,X)$ of coherent $G$-sheaves on a
noetherian $k$-scheme $X$, the group  $K_0(\coh(G,X)^{split})$ is much
too large : in the reductive case, i.e. when $p={\rm char} k \nmid
\#G$ (and in particular when $G=1$), we do not recover the usual
equivariant $K$-theory.

Since the interesting exact structures in the $K$-theory of schemes
comes from abelian categories, we have to reinterpret
the group $K_0(\mc A_{tot}^{split})$ in terms of the Grothendieck
group of an abelian category.

For this purpose, if the ring $\mc A_{tot}$ is right noetherian (see
Definition \ref{ArtNoeth}), the category
$\mod \mc A_{tot}$ is a candidate, since the evaluation
$F \ra F(R)$ provides an exact functor $\mod\mc A_{tot} \ra \mc A_{tot}$
which in turn induces an epimorphism $K_0(\mod\mc A_{tot}) \thra K_0(\mc
A_{tot})$ (this is a consequence of Lemma \ref{RightAdj} applied with $\mc
A'$ the category with only the object $R$). More precisely we have a
commutative diagram :

\xymatrix@R=5pt{
 &&& K_0(\mc A_{tot}^{split})  \ar@{->>}[dr] \ar[dd] &   \\
 &&&                                 &  K_0(\mc A_{tot}) \\
 &&& K_0(\mod \mc A_{tot}) \ar@{->>}[ur] }

\noindent where the two diagonal arrows are the ones already described and the
vertical one is induced by the Yoneda embedding.

That $K_0(\mod\mc A_{tot})$ contains pertinent information is
shown by the following result, which is essentially a
reinterpretation of a theorem of Auslander and Reiten~:

\begin{thm}
\label{Main}
Suppose given a (non necessarily commutative) algebra $R$ over an
algebraically closed field $k$,
with $R$ finite dimensional over $k$
and of finite representation type, and let $\mc A_{tot}=\mod R$.
Then :

(i) the ring $\mc A_{tot}$ is right noetherian

(ii) the Yoneda embedding $\mc A_{tot} \ra \mod\mc A_{tot}$ induces an
isomorphism $$K_0(\mc A_{tot}^{split}) \simeq K_0(\mod\mc A_{tot})$$

(iii) two $R$-modules of finite type $V$,$V'$ are isomorphic
if and only if in $ K_0(\mod\mc A_{tot})$~:
$$[\mc A_{tot} (\cdot, V)]=[\mc A_{tot} (\cdot, V')]$$
\end{thm}

The rest of this paragraph is devoted to a proof of the Theorem,
which is split in the next four sections. Some of the results are
stronger than strictly needed, because we intend to apply them in the
more general context where the algebra $R$ is not of finite
representation type, and $\mc A \subset \mod R$  has a finite set of
generators.

\subsubsection{$K_0(\mc A_{tot}^{split})$ is abelian free of finite rank}

\begin{lem}
\label{KS}
Let $\mc A$ be a full subcategory of $\mc A_{tot}$,
and $V$, $V'$ be two objects of $\mc A$.
Then $[V]=[V']$ in $K_0(\mc A^{split})$ if and only if $V\simeq V'$.
\end{lem}

\begin{proof}
Starting from $[V]=[V']$ we easily get the existence of an object $W$
of $\mc A$ such that $V \oplus W \simeq V' \oplus W$. But then the classical
Krull-Schimdt theorem (see for instance \cite{Lam}, Corollary 19.22)
allows to say that in fact $V\simeq V'$.
\end{proof}

The Lemma \ref{KS} proves the part (iii) of Theorem \ref{Main} is a
consequence of part (ii).
Moreover :

\begin{prop} If $\mc A$ be a projectively complete
  full subcategory of $\mc A_{tot}$, then
  $K_0(\mc A ^{split})$ is isomorphic to the free abelian group generated by
isomorphism classes of indecomposables objects of $\mc A$.
\end{prop}

\begin{proof}
Because $\mc A$ be is projectively complete, the notion of
indecomposable in $\mc A$ is the same as the one in $\mc A_{tot}$, and
we thus can use the Krull-Schimdt theorem again to define a morphism from
$K_0(\mc A ^{split})$ into this free group. The fact that it is an
isomorphism is immediate.
\end{proof}

\subsubsection{$K_0(\mod\mc A_{tot})$ is abelian free}

 We begin by a few definitions concerning an arbitrary additive
category taken from \cite{Mitch}~:

\begin{defi}
\label{ArtNoeth}
An additive category $\mc A$ is called  \emph{right artinian}
(resp. \emph{right noetherian}, resp. \emph{semi-simple})
when for each object $V$ of $\mc A$, the functor $\mc A(\cdot, V)$
is an artinian(resp. noetherian, resp semi-simple)
object of $\Mod \mc A$.
\end{defi}

The dual (left) notion is obtained as usual by replacing $\mc A$ by
$\mc A^{op}$.

Now we recall the definition of the Kelly radical (see \cite{RadKel}).

\begin{defi}
\label{KelRad}
The \emph{Kelly radical} $\rad \mc A$ of an additive category
$\mc A$ is the two-sided ideal of $\mc A$ defined by :
$$\rad \mc A (V,V')=\{ f \in \mc A(V,V')/\; \forall g \in \mc A(V',V)\;\;
1_V -gf \;is\; invertible\}$$
for all pair of objects $(V,V')$ of $\mc A$.
\end{defi}

In \cite{Mitch} it is shown that, as in the one object case, the
notions of left and right semisimplicity coincide, but that however a
right artinian ring need not be right noetherian (the problem being
that $\rad \mc A$ is not necessarily nilpotent). In consequence
B.Mitchell suggests the following definition.

\begin{defi}
An additive category $\mc A$ is said \emph{semi-primary} if

(i) $\mc A/\rad \mc A$ is semi-simple,

(ii) $\rad \mc A$ is nilpotent.

\end{defi}

The advantage of this definition is that the classical
Hopkins-Levitzki Theorem holds now with several objects :

\begin{prop}
\label{HL}
Let $\mc A$ be a semi-primary additive category, and $F$ a right $\mc
A$-module. The following are equivalent :

(i) $F$ is noetherian

(ii) $F$ is artinian

(iii) $F$ is of finite length.

\end{prop}

\begin{proof}
The classical proof applies unchanged : see for instance \cite{Lam} \S
4.15.
\end{proof}

\begin{cor}
\label{ft->fl}
Let $\mc A$ be an additive category. If $\mc A$ is right artinian
and $\rad \mc A$ is nilpotent, then any $\mc A$-module of finite type is of
finite length.
\end{cor}

\begin{proof}
The category $\mc A/\rad \mc A$ is right artinian of zero radical,
hence semi-simple (\cite{Mitch} Theorem 4.4). So $\mc A$ is
semi-primary.

Now if $F$ is a module of finite type, $F$ is artinian, hence of finite
length by Proposition \ref{HL}.
\end{proof}

\begin{lem}
\label{Aus-art}
If $R$ is a finite dimensional $k$-algebra,
and $\mc A$ is subcategory of $\mc A_{tot}=\mod R$
with a finite set of additive generators, then :

(i) the category $\mc A$ is right artinian

(ii) $\rad \mc A$ is nilpotent.

\end{lem}

\begin{proof}
(i)
First note that an arbitrary additive category $\mc A$ is right artinian
if and only if its projective completion $Q\mc A$ is. By hypothesis,
there exists a one object full subcategory $\mc
A'$ of $\mc A_{tot}=\mod R$ such that $Q\mc A' = Q\mc A$,
and $\mc A'$ is obviously right artinian, since it is a
finite dimensional $k$-algebra.

(ii) the restriction along $\mc A \ra Q\mc
A$ provides a bijection between two-sided ideals of $Q\mc A$
and two-sided ideals of $\mc A$ (\cite{Street} Proposition 2). This bijection is
compatible with the product of ideals, and sends $\rad Q\mc A$ to $\rad
\mc A$  (\cite{Street} Proposition 10). With the notations of the
proof of (i), it is thus enough to show that $\rad \mc A'$ is nilpotent,
but this is only the classical fact that the Jacobson radical of an
artinian ring is nilpotent (see \cite{Lam} Theorem 4.12).

\end{proof}

\begin{prop}
\label{simplgen}
If $R$ is a finite dimensional $k$-algebra, and $\mc A$ is subcategory of $\mc A_{tot}=\mod R$
with a finite set of additive generators, then
$\mod \mc A$ is an abelian category, and every object of $\mod \mc
A$ is of finite length. Hence the group $K_0(\mod \mc A)$ is abelian
free, generated by the classes of simple objects in $\mod \mc A$.
\end{prop}

\begin{proof}
Since $\mc A$ is Morita equivalent to an one object right artinian (hence
right noetherian) ring, $\mod \mc A$ is an abelian category.
The second assertion results from Corollary \ref{ft->fl} and Lemma
\ref{Aus-art}. The third assertion is a consequence of the
d\'evissage Theorem (see \cite{Quil} \S 5 Theorem 4 Corollary 1)
for $K_0$, i.e. of the Jordan-H\"older Theorem.
\end{proof}

\subsubsection{${\rm rk}(K_0(\mc A_{tot}^{split}))= {\rm rk}(K_0(\mod\mc A_{tot}))$}

The justification of this equality is the well known remark which
pushed Auslander to introduce functor categories in the study of
representation of Artin algebras : there is a one to one
correspondence between isomorphism classes of
\emph{indecomposables} right $R$-modules of finite type and
isomorphism classes of \emph{simples} right $\mc A_{tot}$-modules
of finite type. This is easily proved directly (see \cite{Gabi} \S
1.2). In this paragraph we give an interpretation of this fact in
a more general context.

According to Proposition \ref{simplgen}, we have, if $R$ is of
finite representation type : $K_0(\mod \mc A_{tot})=K_0((\mod \mc
A_{tot})_{ss})$, where $(\mod \mc A_{tot})_{ss}$ is the full
subcategory of $\mod \mc A_{tot}$ whose objects are the
semi-simple right $\mc A_{tot}$-modules. So what is left to do is
to describe $(\mod \mc A_{tot})_{ss}$. In fact we can in fact
describe $(\mod \mc A)_{ss}$ for a \emph{semi-local} category $\mc
A$, in the following sense :

\begin{defi}
A category $\mc A$ is said \emph{semi-local} if the category
$\mc A/\rad\mc A$ is semi-simple.
\end{defi}

This notion is linked with the notion of radical in the following way.

\begin{defi}
Let $\mc A$ be a category and $F$ a right $\mc A$-module.
We denote by $\rad F$ and call the \emph{radical} of $F$ the
intersection of all maximal submodules of $F$.
\end{defi}

\begin{lem}
\label{rad=rad}
Let $\mc A$ be a semi-local category and $F$ a right $\mc
A$-module. Then $\rad F=F \rad \mc A$.
\end{lem}

\begin{proof}
The usual proof (\cite{Lam} \S 24.4) applies without change.
\end{proof}

\begin{prop}
\label{ss}
Let $\mc A$ be a semi-local category. There is a natural isomorphism~:
$$Q(\mc A/\rad\mc A)\simeq (\mod \mc A)_{ss}$$
\end{prop}

\begin{proof}
Since $\mc A/\rad\mc A$ is semi-simple, the inclusion $Q(\mc A/\rad\mc
A) \subset \mod (\mc A/\rad\mc A)$ is an equality. So the Proposition will
follow from the next Lemma.
\end{proof}

\begin{lem}
Let $\mc A$ be a semi-local category and $\mc B =\mc A/\rad\mc A$.

(i) The functor $R:\Mod \mc B \ra \Mod \mc A$ induced by $\mc A \thra \mc B$ is
fully faithful.

(ii) The image of $\mod \mc B$ under $R$ is $(\mod \mc A)_{ss}$
\end{lem}

\begin{proof}
(i) According to Lemma \ref{rad=rad}, we can define a functor
$K:\Mod \mc A \ra \Mod \mc B$ by setting on objects $K(F)=F/\rad F$. This is a
left adjoint of $R$ (a left Kan extension along $\mc A \thra \mc B$),
and since $KR \simeq 1_{\mc B}$, $R$ is fully faithful (see for
instance \cite{work} Chapter IV \S 3 Theorem 1).

(ii) To see first that $(\mod \mc A)_{ss}\subset R(\mod \mc B) $, let $F$
be an object of $(\mod \mc A)_{ss}$. Then according to Lemma \ref{rad=rad}
and \cite{Street} Proposition 9, we have $\rad F = F \rad \mc A = 0$,
hence $F \simeq F/\rad F=RK(F)$.

To show the opposite inclusion $R(\mod \mc B)\subset(\mod \mc A)_{ss}$,
fix $G$ a semi-simple right $\mc B$-module. Since $R$ is additive we can
assume that $G$ is in fact simple. Let $F\subset R(G)$ be a right $\mc
A$ submodule of finite type.
Then we have a diagram :

\xymatrix@R=5pt{
 &&&  F \ar@{->>}[ddr] \ar@{^{(}->}[rr] &&        G   \\
 &&&&&\\
 &&&                     & F/\rad F=K(F) \ar[uur]&
 }

Since $G$ is simple, and $K(F) \ra G$ is a monomorphism, there are only two
possibilities : either $K(F)=G$ , and then $F=R(G)$, or $K(F)=0$,
which implies, by Lemma \ref{rad=rad}, Nakayama's Lemma (which is
still valid for rings with several objects), and the fact that $F$ is
finitely generated, that $F=0$. Hence $R(G)$ has no proper subobject
in $\mod \mc A$ and by definition is in $(\mod \mc A)_{ss}$.

\end{proof}

Now we can apply Proposition \ref{ss} to the case where $\mc A$ is a
projectively complete subcategory of $\mc A_{tot} = \mod R$, without
finiteness assumption.
To show that $\mc A$ is semi-local, choose an skeleton
$\mc S$ of the category of indecomposables in $\mc A$,
possibly infinite, so that $Q\mc S =\mc A$.

\begin{lem} With notations as above
$$Q(\mc S/\rad \mc S) =\mc A/\rad \mc A$$
\end{lem}

\begin{proof}
Thanks to \cite{Street}, Proposition 10, we can identify  $\mc S/\rad
\mc S$ to a whole subcategory of $\mc A/\rad \mc A$. But now the equality of
the Lemma is obvious, since the functor $\mc A \ra \mc A/\rad \mc A$ is
additive, hence preserves direct sums.
\end{proof}

Now, since semi-simplicity is Morita invariant (a ring is semi-simple
if and only if all its modules are), we are reduced to show the
semi-simplicity of $\mc S/\rad \mc S$.

Recall that a \emph{corpoid} is an additive category where all non zeros
maps are invertible. Given a field $k$ and a set $S$, we have a
associated corpoid, denoted by $kS$, defined by
$kS(V,V)=k$ for any object $V$, and $kS(V,V')=0$ for $V\neq V'$.

\begin{lem}
\label{corpo}
With notations as above, let $S$ be the underlying set of $\mc
S$. Then :
$$\mc S/\rad \mc S = kS $$
\end{lem}

\begin{proof}
For any object $V$ of $\mc S$, the ring $\mc S(V,V)$ is local (see
\cite{Lam}, Theorem 19.17). Hence $\mc S(V,V)/\rad \mc S (V,V)$ is a
skew field, finite dimensional over $k$. Since $k$ is supposed to be
algebraically closed, it must be $k$ itself.

Now consider two different objects $V$,$V'$ of $\mc S$. We have to show
that $\rad \mc S (V,V')=\mc S(V,V')$. Suppose this is not the case,
and let $f$ be an element of $\mc S(V,V')$ be not in $\rad \mc S (V,V')$.
According to the definition of the Kelly radical (see \ref{KelRad}),
there exists $g$ in
$\mc S(V',V)$ such that $1_V-gf$ is not invertible in $\mc S(V,V)$.
Because $\mc S(V,V)$ is local, $gf$ must be invertible (see
\cite{Lam}, Theorem 19.1). Since $V'$ is
indecomposable, this implies that $V \simeq V'$, and this contradicts the
definition of $\mc S$.
\end{proof}

Since $\mc S/\rad \mc S$ has zero radical, and is immediately
seen, thanks to Lemma \ref{corpo}, as right artinian, it is
semi-simple (\cite{Mitch} Theorem 4.4). Hence $\mc A$ is
semi-local, and combining Proposition \ref{ss} and Lemma
\ref{corpo} we get :

\begin{prop}
\label{ind<->simpl}
Let $R$ a finite dimensional $k$-algebra, $\mc A$ a
projectively complete subcategory of $\mc A_{tot} = \mod R$.
Let also $S$ be the set of isomorphism classes of indecomposables
right $R$-modules of finite type contained in  $\mc A$.
Then there is a natural isomorphism :
$$Q(kS) \simeq (\mod \mc A)_{ss}$$
In particular $S$ is in one to one correspondence with the set of
isomorphism classes of simple objects in $\mod \mc A$.

\end{prop}

\subsubsection{The morphism $K_0^{split}(\mc A_{tot}) \ra K_0(\mod\mc A_{tot})$ is
an epimorphism}

\begin{thm}[Auslander-Reiten]
\label{Aus-Rei}
Let $R$ be a finite dimensional $k$-algebra, and $\mc A_{tot} = \mod R$.
Then every simple right $\mc A_{tot}$-module $F$ admits a projective
resolution of the form :
$$0 \ra \mc A_{tot}(\cdot, V'')\ra \mc A_{tot}(\cdot, V')\ra \mc A_{tot}(\cdot, V)\ra
F \ra 0$$
\end{thm}

\begin{proof}
See \cite{Gabi}, \S 1.3.
\end{proof}

\begin{cor}
Let $R$ be a finite dimensional $k$-algebra of finite type. Then
the morphism $K_0^{split}(\mc A_{tot}) \ra K_0(\mod\mc A_{tot})$ is
an epimorphism.
\end{cor}

\begin{proof}
This is a direct consequence of Proposition \ref{simplgen} and of
Theorem \ref{Aus-Rei}.
\end{proof}

\subsubsection{k-categories and additive categories}

Since we started with a $k$-algebra $R$, we could have worked with
$k$-categories, i.e. categories enriched in $\Mod k$, rather
than with additive categories.

However, for modules categories, it
does not make a significant difference. Indeed, for a $k$-category $\mc A$,
denote by $\mc A_0$ the underlying additive category. Then is it easy
to see that restriction along the forgetful functor $\Mod k\ra
\ab$ induces an equivalence~:

$$[\mc A^{op}, \Mod k]_0 \simeq [\mc A^{op}_0, \ab]$$

So we can enrich our functor categories to see them as $k$-categories.
This formulation has some advantages, in particular :

\begin{lem}
\label{CohCoh}
Let $\mc A$ a be category enriched over $\mod k$, admitting a finite
set of additive generators. Then :

$$[\mc A^{op}, \mod k]_0 \simeq \mod (\mc A_0)$$
\end{lem}

\begin{proof}
This is a direct consequence of the Yoneda Lemma.
\end{proof}

This applies in particular for $\mc A_{tot}=\mod R$, for $R$ a finite
dimensional $k$-algebra of finite representation type.

\subsection{The general case}

\subsubsection{What is left}

Let $R$ be a finite dimensional $k$-algebra. If one makes no
assumption on the representation type of $R$, one can not really work
with the ring $\mc A_{tot}$ any longer, since it needs not be right
noetherian. Instead we fix a projectively complete subcategory $\mc A$
of $\mc A_{tot}$, admitting a finite set of additive generators. For
such a category $\mc A$, Propositions \ref{simplgen} and
\ref{ind<->simpl} hold, which shows that $K_0(\mod \mc A)$ and
$K_0(\mc A^{split})$ are free abelian groups of the same finite rank.
However, we can not say that the morphism induced by Yoneda between
these two groups is an isomorphism.

We can be more precise : if we suppose that $\mc A$ contains the free
object $R$, hence the category of projective right $R$-modules $\mc P
=QR$, then we have a natural commutative diagram :

\xymatrix{
 &&& K_0(\mc A^{split})  \ar[r] \ar[d] &  K_0(\mc
 A_{tot}^{split})\ar@{->>}[d]\\
 &&&        K_0(\mod \mc A)     \ar@{->>}[r]          &  K_0(\mc
 A_{tot})
}

\subsubsection{Inverse of the devissage isomorphism}

It is convenient to give an explicit description of the inverse isomorphism
of the one given by devissage, so we introduce some notations.

Again we fix a projectively complete subcategory $\mc A$
of $\mc A_{tot}=\mod R$, admitting a finite set of additive
generators, $\mc S$ a skeleton of the subcategory
of indecomposables of $\mc A$, and $S$ the underlying finite set.

\begin{defi}
For each object $I$ of $\mc S$, we denote by $S_I$ the simple object
of $\mod \mc A$ consisting of the quotient of the projective functor
$\udl I:=\mc A(\cdot, I)$ by the intersection ${\rm rad}\, I$ of
its maximal subobjects.
\end{defi}

\begin{lem}
\label{expl}
Let $F$ be an object of $\mod \mc A$.

We have in
$K_0(\mod \mc A)$ :

$$[F]=\sum_{I\in S}\dim_k F(I)[S_I]$$
\end{lem}

\begin{proof}
We have a canonical morphism $K_0(\mod \mc A)\ra {\rm Map}( S,
\mathbb Z)$ sending $[F]$ to $I \ra \dim_k F(I)$. This map sends $[S_J]$
to $$I \ra \dim_k \mod \mc A (\udl{I},S_J)=I \ra \dim_k \mod \mc A
(S_I,S_J)= I \ra \delta_{I,J}$$
hence it is in fact an isomorphism, and the formula follows.
\end{proof}


\section{Enriched Kan extension}
\label{Kan}

We give a brief summary of the notions of enriched category theory
we will need in the sequel. We follow essentially the lines of the
foundational papers \cite{EilKel}, \cite{Kelly}, \cite{DayKel},
by B.Day, G.M. Kelly and S.Eilenberg.

\subsection{Notations}

As in \cite{Kelly}, \cite{DayKel}, a \emph{closed category}
$\mc S$ will denote more precisely a symmetric closed monoidal
category as defined in \cite{EilKel}.
The definitions of $\mc S$-category, $\mc S$-functor,
$\mc S$-natural transformation used are also those given in \cite{EilKel}.
We will often stress the presence of a structure of
closed or enriched category by using bold letters.

We will write $\mc S-CAT$ for the hypercategory (or strict
$2$-category) whose objects are the $\mc S$-categories, $1$-arrows are
the $\mc S$-functors, and $2$-arrows are the $\mc S$-natural
transformations.
In particular we define as usual
${\bf \mc C}at = {\bf \mc E}ns-CAT$ and ${\bf \mc H}yp = {\bf \mc C}at-CAT$.

We will write ${\bf \mc C}l$ for the $2$-category of closed categories.

The notion of adjunction between $\mc S$-functors used is the one defined
in \cite{Kelly}, \S2. It implies the existence of an adjunction
between the underlying functors, but is not implied by this one.

Let $\mc B$ be a $\mc S$-category, $X \in {\rm Obj} \mc S$, $B \in
{\rm Obj} \mc B$.

The tensor $X\otimes B \in {\rm Obj} \mc B$ is
characterized by the existence of a $\mc S$-natural isomorphism :

$$\mc B( X\otimes B, B') \simeq [X, \mc B (B,B')]$$
where the brackets in the right-hand side denote the internal Hom of
$\mc S$.

Dually the cotensor $[X,B] \in {\rm Obj} \mc B$ is
characterized by the existence of a $\mc S$-natural isomorphism :

$$\mc B(B',[X,B]) \simeq [X, \mc B (B',B)]$$
(see \cite{Kelly}, \S4).

The definition of a complete (resp. cocomplete)
$\mc S$-category $\mc B$ is the one given in
\cite{DayKel}, \S2 (this should of course not be confused with the
notion of projective completion of an additive category introduced
above). It implies in particular the existence of cotensor (resp. tensor)
objects, and the existence of small ends (resp. coends), see
\cite{DayKel}, \S 3.3.

Finally, if $\mc B$ and $\mc B'$ are $\mc S$-categories, with $\mc B'$
small and $\mc S$ complete, the category of
$\mc S$-functors between $\mc B'$ and $\mc B$,
and $\mc S$-natural transformations between these functors, can be
enriched in a $\mc S$-category $[\mc B',\mc B]$ by setting, for two
$\mc S$-functors $T,S$~:

$$[\mc B',\mc B](T,S) =\int_{B'} \mc B (TB', SB')$$
(see \cite{DayKel}, \S4).

If $\mc B$ is complete (resp. cocomplete) then $[\mc B',\mc B]$ is
also complete (resp. cocomplete), and limits (resp. colimits) are
formed termwise (see \cite{Basic}, Chapter 3, \S 3.3).

\subsection{Enriched left Kan extension}

\begin{prop}[Day-Kelly]
\label{DayKellyleft}
Let $\mc S$ be a closed category, $i: \mc B' \ra \mc B$ a $\mc
S$-functor, where $\mc B$ and $\mc B'$ are small, and $\mc C$ a
cocomplete $\mc S$-category. Then the $\mc S$-functor $R=[i,1]:[\mc B, \mc C]
\ra [\mc B', \mc C]$ admits as a left adjoint the functor $Q$
given on objects by, for $F \in {\rm Obj}[\mc B', \mc C]$ :

$$Q(F)=\int^{B'}\mc B(\cdot,iB')\otimes FB'$$
\end{prop}

\begin{proof}
See \cite{DayKel}, \S 6.1.
\end{proof}

\subsection{Enriched right Kan extension}

\begin{prop}
\label{DayKelly}
Let $\mc S$ be a closed category, $i: \mc B' \ra \mc B$ a $\mc
S$-functor, where $\mc B$ and $\mc B'$ are small, and $\mc C$ a
complete $\mc S$-category. Then the $\mc S$-functor $R=[i,1]:[\mc B, \mc C]
\ra [\mc B', \mc C]$ admits as a right adjoint the functor $K$
given on objects by, for $F \in {\rm Obj}[\mc B', \mc C]$ :

$$K(F)=\int_{B'}[\mc B(\cdot,iB'),FB']$$
If moreover $i$ is fully faithful,
the counit of this adjunction is an isomorphism :
$RK \simeq 1$ (equivalently, $K$ is fully faithful).
\end{prop}

\begin{proof}
The first statement is the dual of Proposition
\ref{DayKellyleft}, and the second follows from the enriched Yoneda
Lemma (see \cite{DayKel}, \S5).
\end{proof}

\begin{cor}
\label{EnrRightKan}
Suppose moreover that $\mc C = \mc S$.
Then the $\mc S$-functor $R=[i,1]:[\mc B, \mc S]
\ra [\mc B', \mc S]$ admits as a right adjoint the functor :
\begin{equation*}
 \xymatrix @R=0pt {
  K:[\mc B', \mc S] \ar[r]  &  [\mc B, \mc S]& \\
          F        \ar[r] &  (V \ar[r] & [\mc B', \mc
  S](\mc B (V,i\cdot),F))}
\end{equation*}
\end{cor}

\begin{proof}
By definition $K(F)(V)=\int_{B'}[\mc B(V,iB'),FB']$. Since $\mc S$
is symmetric, the cotensor $[\cdot, \cdot]$ in $\mc S$ coincide with
the internal Hom of $\mc S$. The formula given follows from the definition of
the $\mc S$-category $[\mc B', \mc S]$.
\end{proof}

\begin{cor}
\label{EnrRightKan2} Suppose moreover that $\mc C = \mc S$, $i$ is
fully faithful, and every object of $\mc B$ is a retract of an
object of $\mc B'$. Then the adjunction of Corollary
\ref{EnrRightKan} is an equivalence.
\end{cor}

\begin{proof}
Consider the composite embedding $\mc B'^{op} \ra \mc B^{op} \ra
[\mc B, \mc S]$. Since $\mc B^{op}$ is dense in $[\mc B, \mc S]$, it
follows from \cite{Basic}, Chapter 5, Proposition 5.20, that $\mc
B'^{op}$ is dense in $[\mc B, \mc S]$. But \cite{Basic}, Chapter 5,
Theorem 5.1 (ii) and the enriched Yoneda Lemma again show that the
restriction $R:[\mc B, \mc S]\ra [\mc B', \mc S]$ is fully
faithful. It follows now from \cite{Basic}, Chapter 1, \S 1.11 that
the unit $1 \Longrightarrow RK$ of the adjunction is an isomorphism.
\end{proof}

In fact, the part of classical Morita theory briefly described in \S
\ref{Morita} lifts to the general enriched context (see in particular
\cite{Basic}, Chapter 5, Proposition 5.28). Since we will not use these
results, we do not recall them.


\section{Rings with several objects on a scheme}

\label{RingedSchemes}

All schemes considered in the sequel are supposed noetherian.

\subsection{Category of quasi-coherent sheaves on a scheme}

\begin{prop}
Let $Y$ be a scheme, and $\qcoh Y$ the category of quasi-coherent
sheaves on $Y$. There is a closed category $\Qcoh Y$ whose underlying
category is $\qcoh Y$.
\end{prop}

\begin{proof}

We have to give first the seven data (in fact six independent)
defining a closed category as given in \cite{EilKel}, Chapter I,
\S 2, and we choose the natural ones. Checking the axioms CC1 to
CC5 can be done locally, hence deduced from the corresponding
facts for categories of modules over a commutative ring, or even
directly.

Moreover, there is a tensor product on $\qcoh Y$ defined by the existence of
a natural isomorphism
\begin{equation}
\label{tensor}
\qcoh Y(\mc F \otimes \mc G, \mc H) \simeq
\qcoh Y(\mc F, \mc{H}om (\mc G, \mc H))
\end{equation}

\noindent It lifts to a $\qcoh Y$-natural transformation :
$$\mc{H}om (\mc F \otimes \mc G, \mc H) \simeq
\mc{H}om (\mc F, \mc{H}om (\mc G, \mc H))$$
Hence we deduce from \cite{EilKel}, Chapter II, Theorem 5.3, that the closed
category defined above is in fact monoidal.

At least, using the definition \ref{tensor} of the tensor product and
the natural symmetry of $\mc{H}om$

$$\qcoh Y(\mc F, \mc{H}om (\mc G, \mc H)) \simeq
\qcoh Y(\mc G, \mc{H}om (\mc F, \mc H))$$
one defines a symmetry for the tensor product and check the axioms
MC6-MC7 of \cite{EilKel}, \S Chapter III, \S1.

\end{proof}

\begin{prop}
\label{QcohCompl}
$\Qcoh Y$ is complete and cocomplete.
\end{prop}

\begin{proof}
According to the definition given in \cite{DayKel}, \S2, it suffices
to show that $\qcoh Y$ is complete and cocomplete. Since (co)limits
must commute with localization, one first construct them locally,
which is possible, because the categories of modules over a ring are
complete and cocomplete (see \cite{Schub} 7.4.3, 8.4.3).
Since limits are universal, one can glue the local limits together to
get global limits.
\end{proof}

\begin{prop}
\label{quasicohfunc_*}
Let  $f:Y' \ra Y$ be a morphism of schemes.
The functor $f_* : \qcoh Y' \ra \qcoh Y$ can be lifted
in a closed functor $f_* : \Qcoh Y' \ra \Qcoh Y$.



\end{prop}

\begin{proof}
  To lift $f_* : \qcoh Y' \ra \qcoh Y$ one needs to specify a natural
transformation $\hat{f_*}: f_* \mc H om (\mc F',\mc G') \ra
\mc H om(f_*\mc F', f_*\mc G')$, and a morphism $f_*^0 : \mc O_Y \ra f_*
  \mc O _{Y'}$, and once again one chooses the obvious ones.
To check the axioms CF1-CF3 (resp. MF4) given in \cite{EilKel},
Chapter I, \S 3 (resp. Chapter III, \S 1) is a long but easy task.


\end{proof}

\begin{prop}
\label{quasicohfunc^*}
Let  $f: Y' \ra Y$ be a morphism of schemes.
The functor $f^* : \qcoh Y \ra \qcoh Y'$ can be lifted
in a closed functor $f^* : \Qcoh Y \ra \Qcoh Y'$.
\end{prop}

\begin{proof}

To lift $f^* : \qcoh Y' \ra \qcoh Y$, one needs in a similar way
to specify a natural transformation $\hat{f^*}: f^* \mc H om (\mc F,\mc G) \ra
\mc H om(f^*\mc F, f^*\mc G)$, and a morphism ${f^*}^0 : \mc O_Y \ra f^*
  \mc O _{Y'}$. For ${f^*}^0$ one makes the natural choice of the
  inverse of the isomorphism given from $f_*^0$ by adjunction. To
  construct $\hat{f^*}$, one first notices that the composition
$$ f_* \mc Hom (f^*\mc F, \mc G') \ra \mc Hom (f_*f^*\mc F, f_* \mc
  G') \ra \mc Hom(\mc F, f_*\mc G') $$
given by $\hat{f_*}$ and adjunction, is in fact a natural isomorphism.
  One gets then immediately a natural $\hat{f^*}$, also given by
  adjunction. Using the Proposition \ref{quasicohfunc_*},
  one checks the axioms CF1-CF3 and MF4 again.


\end{proof}

\begin{prop}

\label{AdjCl}
Let  $f:Y' \ra Y$ be a morphism of schemes.
The usual adjunction between $f^*$ and $f_*$ in ${\bf \mc C}at$
can be lifted in an adjunction in the $2$-category of closed
categories ${\bf \mc C}l$.
\end{prop}

\begin{proof}
One checks that the unit and counit of the adjunction are in fact
closed natural transformations, i.e. that they verify the axioms CN1
and CN2 given in \cite{EilKel}, Chapter I, \S 4.
\end{proof}

\begin{cor}
\label{AdjHyp}
The functors $f^*$ and $f_*$ belong to a natural ${\bf \mc
  C}at$-enriched adjunction between $\Qcoh Y-CAT$ and $\Qcoh Y'-CAT$.
\end{cor}

\begin{proof}
According to \cite{Kelly}, \S2,
we can just push the adjunction of Proposition \ref{AdjCl} along the
canonical $2$-functor ${\bf \mc C}l \ra {\bf \mc H}yp$ sending $\mc S$
to $\mc S-CAT$, to get an adjunction in ${\bf \mc H}yp$, i.e., by
definition, a ${\bf \mc C}at$-enriched adjunction.
\end{proof}

\begin{cor}

 \label{CorAdjCl}
There is a natural $2$-arrow in the $2$-category
${\bf \mc C}l$ of closed categories~:

  \xymatrix@R=12pt{
    &&&&& \Qcoh Y' \ar[rrd]^{H^0(Y',.)}  \\
    &&&&&               && \ab \\
    &&&&& \Qcoh Y \ar[uu]^{f^*} \ar[rru]_{H^0(Y,.)}="A"  \ar@{<=}[uu];"A"}

\end{cor}

\begin{proof}
This is obtained by pushing in the $2$-category
  $\mc Cl$ the unit of the adjunction of Proposition \ref{AdjCl}
along the closed functor $H^0(Y,.)$.
\end{proof}





\subsection{Ringed schemes}

\begin{defi}

Let $Y$ be a scheme.

(i) A \emph{ring (with several objects) on $Y$} is, by definition, a category
$\mc A$ enriched on $\Qcoh Y$.
The pair $(Y, \mc A)$ is  called a \emph{ringed scheme}.

(ii) A morphism of ringed schemes $(Y', \mc A')\ra (Y, \mc A)$ is a couple
$(f, f^\#)$ where
$f : Y' \ra Y$ is a scheme morphism and $f^\# : \mc A \ra f_*\mc A'$
is a morphism of $\Qcoh Y$-categories.

\end{defi}

\begin{rem}
There is an obvious notion of $2$-arrow
between two ringed schemes morphisms of same source and target, and
ringed schemes thus form a $2$-category.
\end{rem}

A natural operation to consider is, for every open $i: U \ra Y$,
to push $\mc A$ via $H^0(U,\cdot) \circ i^* :
\Qcoh Y \ra \Qcoh U \ra \ab$, which gives an
additive category $\mc A(U)$. Moreover, if $i': U' \ra Y$ is another
open, every $Y$-inclusion $f:U' \ra U$
gives an additive restriction functor $\mc A(U) \ra \mc
A(U')$. Indeed from the lax-functor $(Sch/Y)^{op} \ra {\bf \mc C}l$ sending
$i: U \ra Y$ to $\Qcoh U$, and Corollary \ref{CorAdjCl}, we get the
diagram in ${\bf \mc C}l$ :

\xymatrix@R=12pt{
    &&&&& \Qcoh U' \ar[rrd]^{H^0(U',.)}  \\
    &&&\Qcoh Y \ar[rru]^{i'^*}="B" \ar@{=>}[rrd];"B"
\ar[rrd]_{i^*} &&               && \ab \\
    &&&&& \Qcoh U \ar[uu]^{f^*} \ar[rru]_{H^0(U,.)}="A"  \ar@{<=}[uu];"A"}

We can then push the resulting $2$-arrow via the canonical $2$-functor
${\bf \mc C}l \ra {\bf \mc H}yp$ sending $\mc S$ to $\mc S-CAT$, and
  then evaluate at $\mc A$. In this way, $\mc A$ can be seen as a sheaf
  of rings with several objects on $Y$.







\subsection{Quasicoherent sheaves on ringed schemes}

\subsubsection{Definition}

\begin{defi}
Let $(Y,\mc A)$ be a ringed scheme.

(i) A quasicoherent sheaf on $Y$ is by definition an enriched functor
from $\mc A^{op}$ to $\Qcoh Y$.

(ii) A morphism of $\mc A$-sheaves is an enriched
natural transformation.

(iii) We denote by $\Qcoh(Y,\mc A)=[\mc A^{op}, \Qcoh Y]$ the
corresponding enriched category, and by $\qcoh(Y,\mc A)$ the underlying category.

\end{defi}

When $Y={\rm Spec}\,B$ is affine, then clearly the global sections
functor induces an equivalence $\Qcoh(Y,\mc A)\simeq  [\mc
A^{op}(Y),B \,\Mod]$ as $\Qcoh Y \simeq B \,\Mod$ enriched
categories.

\subsubsection{Functoriality : definition}

\label{ChangeRing}

To define push-forward and pull-back for sheaves on ringed schemes
we need some notations.

First consider the simplest case of a morphism $(1, i):
(Y,\mc A')\ra(Y,\mc A)$.
Because of
Proposition \ref{QcohCompl} we can apply the results quoted in part
\ref{Kan}. Thus we get left and right adjoints for the restriction
functor  $\Qcoh(Y,\mc A) \ra \Qcoh(Y,\mc A')$. The left adjoint
will be denoted as usual by $\otimes_ {\mc A'} \mc A$.

Now consider a general morphism $(f, f^\#):(Y',\mc A')\ra(Y,\mc A)$.

We denote by $adj$ the ${\bf \mc C}at$-natural
isomorphism  $$adj :
\Qcoh Y-CAT(\mc B, f_* \mc B') \simeq \Qcoh Y'-CAT(f^*\mc B, \mc B') $$
given by Corollary \ref{AdjHyp}. We write $\epsilon$ (resp. $\mu$)
for the counit (resp. for the unit) of the adjunction between $f^*$
and $f_*$ in ${\bf \mc C}l$.

Recall that if $\Phi : \mc S' \ra \mc S$ is any closed functor,
there is an associated $\mc S$-functor $\Phi\mc S' \ra \mc S$, and
we will denote it by $cr(\Phi)$ (see \cite{EilKel}, Chapter I, Theorem 6.6).

\begin{defi}
\label{FunctRingedSchemes}
Let  $(f, f^\#):(Y',\mc A')\ra(Y,\mc A)$ be a morphism of ringed schemes.

(i) We define $f_{\triangle} : f_*\Qcoh(Y',\mc A')\ra \Qcoh(Y,\mc A)$ as
the $\Qcoh Y$-functor making the following diagram commute :

\begin{equation*}
\xymatrix {
   f_*[\mc A'^{op},\Qcoh Y'] \ar[r]^{f_*}
   \ar[ddr]_{f_{\triangle}} &
   [f_*\mc A'^{op},f_*\Qcoh Y']
     \ar[dd]^{[f^{\#},cr(f_*)]}\\
   &                            \\
    & [\mc A^{op}, \Qcoh Y]
   }
\end{equation*}

(ii) We define $f^{\triangle} : \Qcoh(Y,\mc A)\ra
f_*\Qcoh(Y',\mc A')$ as the $\Qcoh Y$-functor making the following diagram commute :
\begin{equation*}
\xymatrix {
   [\mc A'^{op},\Qcoh Y'] & \\
   &                            \\
 [f^*\mc A^{op}, \Qcoh Y']
 \ar[uu]^{\otimes_{f^*\mc A}\mc
   A'} &         \\
   &                            \\
   [f^*\mc A^{op}, f^*\Qcoh Y]
   \ar[uu]^{[1,cr(f^*)]}
   &  f^*[\mc A^{op}, \Qcoh Y] \ar[l]^{f^*}
   \ar[uuuul]_{adj(f^{\triangle})}
   }
\end{equation*}

\end{defi}

\subsubsection{Functoriality : adjunction}

\begin{prop}
\label{AdjunctionRingedSchemes}
Let  $(f, f^\#):(Y',\mc A')\ra(Y,\mc A)$ be a morphism of ringed schemes.
The couple $(f^{\triangle},f_\triangle)$ is part of a
$\Qcoh Y$-adjunction between $\Qcoh(Y,\mc A)$ and $f_*\Qcoh(Y',\mc A')$.
\end{prop}

\begin{proof}
This is a consequence of the five following lemmas.

\begin{lem}
The ${\bf \mc C}at$-natural
isomorphism  $$adj :
\Qcoh Y-CAT(\mc B, f_* \mc B') \simeq \Qcoh Y'-CAT(f^*\mc B, \mc B') $$
given by Corollary \ref{AdjHyp} lifts to a $\Qcoh Y-CAT$-natural isomorphism
$$\widehat{adj}:[\mc B, f_* \mc B']\simeq f_*[f^*\mc B, \mc B']$$
\end{lem}

\begin{proof}
There are natural morphisms :

\xymatrix{
  f_*[f^*\mc B,\mc B'] \ar[r]^{f_*} & [f_*f^*\mc B,f_* \mc B']
  \ar[r]^{[\mu_{\mc B},1]} &  [\mc B,f_* \mc B'] }

and

\xymatrix{
  f^*[\mc B,f_*\mc B'] \ar[r]^{f^*} & [f^*\mc B,f^*f_* \mc B']
  \ar[r]^{[1,\epsilon_{\mc B'}]} &  [f^*\mc B, \mc B'] }

\noindent and the morphism associated to the second one by adjunction is an
inverse of the first one.
\end{proof}

\begin{lem}

  \label{BigLemma}
 The following diagram in $\Qcoh Y-CAT$ is commutative :

 \begin{equation*}
 \newcommand{\eq}[1][r]
   {\ar@<-3pt>@{-}[#1]
    \ar@<-1pt>@{}[#1]|<{}="gauche"
    \ar@<+0pt>@{}[#1]|-{}="milieu"
    \ar@<+1pt>@{}[#1]|>{}="droite"
    \ar@/^2pt/@{-}"gauche";"milieu"
    \ar@/_2pt/@{-}"milieu";"droite"}
 \xymatrix {
   f_*[\mc A'^{op},\Qcoh Y']
   \ar[dd]_{f_*{[adj(f^\#),1]}} \ar[r]^{f_*}
    &
   [f_*\mc A'^{op},f_*\Qcoh Y']
     \ar[dd]^{[{f}^{\#},1]} \\
   &                            \\
    f_* [f^*\mc A^{op}, \Qcoh Y']
  \eq[r]^{\widehat{adj}}&
   [\mc A^{op},f_*\Qcoh Y']\\
   &                            \\
   f_*[f^*\mc A^{op}, f^*\Qcoh Y]
   \ar[uu]^{f_*{[1, cr({f}^*)]}}
   &  [\mc A^{op}, \Qcoh Y] \ar[l]^{adj^{-1}(f^*)}
   \ar[uu]_
      {[1, adj^{-1}(cr({f}^*))]}
   }
 \end{equation*}
\end{lem}

\begin{proof}

Let us consider first the following subdivision of the top square.

\begin{equation*}
\xymatrix { f_*[\mc A'^{op},\Qcoh Y']
   \ar[dd]_{f_*{[adj(f^\#),1]}} \ar[rr]^{f_*}
    &&
   [f_*\mc A'^{op},f_*\Qcoh Y']\ar@{-->}[ld]_{[f_* adj(f^\#),1]}
     \ar[dd]^{[{f}^{\#},1]} \\
   &     [f_* f^*\mc A^{op},f_* \Qcoh Y'] \ar@{-->}[rd]^{[\mu_{\mc A},1]}
     &                    \\
    f_* [f^*\mc A^{op}, \Qcoh Y'] \ar@{-->}[ru]^{f_*}
  \ar[rr]_{\widehat{adj}^{-1}}&&
   [\mc A^{op},f_*\Qcoh Y']
   }
\end{equation*}
The quadrilateral commutes because of the $\Qcoh Y-CAT$-naturality of $f_*$.
 The bottom triangle commutes by definition of $\widehat{adj}$, and the
 right triangle by definition of $adj$.

 So the top square commutes, and the bottom square is treated in a
 similar way.
\end{proof}

\begin{lem}
The couple $(\otimes_{{f}^*\mc A}\mc
   A',[adj({f}^\#),1]))$ is part of a natural adjunction.
\end{lem}

\begin{proof}
 This is by definition of the tensor product,
 as given in \S \ref{ChangeRing}.
\end{proof}

Note that we can push this adjunction along the $2$-functor
$f_* : \Qcoh Y'-CAT \ra \Qcoh Y-CAT$, to get a first adjunction in the
diagram of the Lemma \ref{BigLemma}. The following Lemma gives a
second one.

\begin{lem}
  The couple
  $([1, adj^{-1}(cr({f}^*))] ,[1,cr({f}_*)])$
  is part of a natural adjunction.
\end{lem}

\begin{proof}
It is of course enough to show that $(adj^{-1}(cr({f}^*)
  ,cr({f}_*))$ is a part of a natural adjunction
  in $\Qcoh Y -CAT$, since we can
  then push it along the $2$-functor
  $[\mc A^{op}, \cdot]:\Qcoh Y -CAT
  \ra \Qcoh Y -CAT$. Hence the following Lemma allows to conclude :
\end{proof}

\begin{lem}
Let $\mc Q$ and $\mc Q'$ be two closed categories,
$L:\mc Q \ra \mc Q'$ and $R:\mc Q' \ra \mc Q$ two
closed functors part of closed adjunction $(L,R,\epsilon,\mu)$
(where $\epsilon$ is the counit and $\mu$ the unit), and $adj : \mc Q-CAT
(\mc B,R\mc B') \simeq \mc Q'-CAT (L\mc B, \mc B')$ the induced ${\bf
  \mc C}at$-natural isomorphism. Then the pair
$(adj^{-1}(cr(L)),cr(R))$ is part of a natural $\mc Q$-adjunction.
  \end{lem}

\begin{proof}
Applying \cite{EilKel}, Chapter I, Proposition 8.10 to the unit :
$\mu : 1 \Longrightarrow RL : \mc Q \ra \mc Q $, we get a $2$-arrow in
$\mc Q-CAT$ :
$\mu : 1 \Longrightarrow cr(RL)\mu_{\mc Q} : \mc Q \ra \mc Q $.
But since the following diagram

\begin{equation*}
\xymatrix{
\mc Q \ar[rr]^{\mu_{\mc Q}}\ar[dd]_{adj^{-1}(cr(L))} && RL \mc Q
\ar[dd]^{cr(RL)} \ar@{-->}[ddll]^{R(cr(L))}
\\&&
\\
R\mc Q \ar[rr]_{cr(R)} && \mc Q}
\end{equation*}
commutes (as one sees thanks to the dotted arrow), we get in fact a
(candidate) unit $\mu : 1 \Longrightarrow cr(R) adj^{-1}(cr(L))
: \mc Q \ra \mc Q $. The (candidate) counit is built by the same
argument, using moreover
the ${\bf \mc C}at$-adjunction between $\mc Q-CAT$ and
$\mc Q'-CAT$.
\end{proof}

We can then compose the two adjunctions we have described in diagram of
the Lemma \ref{BigLemma}, and this shows the Proposition.
\end{proof}

\subsubsection{Stalks}

\label{StalksRinged}
Starting from a ringed scheme $(Y, \mc A)$, we have for each point $Q$
of $Y$ a localization $\mc A_Q$ given by $\mc A_Q=i^*\mc A$, where $i:
{\rm spec}\, \mc O_{Y,Q} \ra Y$ is the canonical morphism.
Pulling back along the ringed scheme morphism
$({\rm spec}\, \mc O_{Y,Q},\mc A_Q) \ra (Y, \mc A)$ gives a stalk functor :

$$\xymatrix@R=0pt{
  &&& \Qcoh(Y,\mc A) \ar[r] & \Qcoh({\rm spec}\, \mc O_{Y,Q},\mc A_Q) \\
  &&& \mc F \ar[r] & \mc F_Q }$$

\subsubsection{$\qcoh(Y,\mc A)$ is abelian}

\begin{defi}
Let $(Y,\mc A)$ be a ringed scheme.
For each object $V$ of $\mc A$, we denote by $\langle V \rangle$ the
full subcategory of $\mc A$ containing only the object $V$.
We will write $p_V$ or $\cdot(V)$ for the canonical projection :
$ \Qcoh (Y,\mc A) \ra \Qcoh (Y, \langle V \rangle)$ .
\end{defi}

For each $V$, the category $\qcoh (Y, \langle V \rangle)$ is the category of
quasicoherent sheaves of modules on a (one object !) algebra on $Y$,
hence it is an abelian category.

\begin{prop}
\label{abelian}
(i) The category $\qcoh(Y,\mc A)$ is an abelian category.

(ii) A sequence $\mc F' \ra \mc F \ra \mc F''$ of $\mc A$-sheaves on $Y$ is
exact if and only if for each object $V$ of $\mc A$ the sequence
$\mc F'(V) \ra \mc F(V) \ra \mc F''(V)$ is exact in $\qcoh
(Y, \langle V \rangle)$.
\end{prop}

\begin{proof}

Because of Proposition \ref{QcohCompl}, $\qcoh(Y,\mc A)$ is complete
and cocomplete, and limits are formed termwise. So (ii) follows from (i).
To show (i),
the only thing to prove is that every monomorphism is a kernel, and
the dual assertion.

Consider the functor given by the canonical projections :

$$ p : \qcoh (Y,\mc A) \ra \prod_V\qcoh (Y,\langle V \rangle)$$

It is faithful, and we know from part \ref{ChangeRing}
that each $p_V$ admits a right adjoint,
hence preserve colimits, hence $p$ itself preserve colimits.

Dually, one sees that $p$ also preserve limits.

Now if $f$ is a monomorphism in $\qcoh(Y,\mc A)$, then $p(f)$ is a
monomorphism in $\prod_V\qcoh (Y,\langle V \rangle)$, which is abelian, hence
$p(f)= {\rm Ker}({\rm Coker}\, p(f))=p({\rm Ker}({\rm Coker} f))$, and
because $p$ is faithful $f= {\rm Ker}({\rm Coker} f)$, hence $f$ is a
kernel. The dual assertion follows similarly.

\end{proof}

Since localization commutes with $\cdot(V)$, we get immediately
from
 Proposition \ref{abelian}~:

\begin{prop}
\label{ExactStalks}
The sequence of $\mc A$-sheaves :
$\mc F' \ra \mc F \ra \mc F''$ is exact if and only if for each $Q$
in $Y$ the sequence of stalks
$\mc F'_Q \ra \mc F_Q \ra \mc F''_Q$
is exact.
\end{prop}


\section{Sheaves of modules for an Auslander algebra on a $G$-scheme}

\label{AuslanderAlgebra}

\subsection{Auslander algebras associated to a $G$-scheme over a field}

\subsubsection{Definition}

In the sequel, we fix an algebraically closed field $k$, and a scheme $X$ over
$k$, endowed with an \emph{admissible} action of a finite group $G$,
so that the quotient scheme $Y=X/G$ exists. The quotient morphism will
be denoted by $\pi : X \ra Y=X/G$. We call this data \emph{a
$G$-scheme over $k$}.

We recall briefly the definition of a $G$-sheaf :

\begin{defi}
Let $\mathcal F$ be a quasicoherent sheaf on the $G$-scheme $X$.
A \emph{$G$-linearization} of $\mathcal F$
is the data of a collection $(\psi_g)_{g \in G}$ of sheaf morphisms
$\psi_g : g_* \mathcal F \rightarrow \mathcal
F$ checking the following conditions :

\begin{enumerate}
\item   $\psi_1=1$
\item   $\psi_{hg}=\psi_h \circ h_*(\psi_g)$
  in other words, the following diagram commute :
\end{enumerate}
\xymatrix{
&&&&h_* g_*\mc F \ar@{=}[d]\ar[r]^{h_*\psi_g} & h_* \mc F \ar[r]^{\psi_h}
&\mc F \\
&&&& (hg)_*\mc F \ar[urr]_{\psi_{hg}}           }
A \emph{$G$-sheaf on $X$} is by definition a quasicoherent sheaf on $X$
endowed with a $G$-linearization.
\end{defi}

Morphism of $G$-sheaves are morphism of sheaves commuting with the
action.
In this way, $G$-sheaves on $X$ form an abelian category $\qcoh
(G,X)$.

If $\mc G$ is a sheaf on the quotient $Y$, the sheaf $\pi^* \mc G$ has
a natural structure of $G$-sheaf, and we get a functor $\pi^* :
\qcoh Y \ra \qcoh(G,X)$. This functor admits as usual a right adjoint,
denoted by $\pi_*^G$, which to a $G$-sheaf $\mc F$ associates the sheaf
on $Y$ given by $U\ \ra (\mc F (\pi^{-1}U))^G$.

\begin{lem}
  \label{QcohGX}
Let $X$ be a $G$-scheme.
The closed structure on $\Qcoh X$ induce on $\qcoh(G,X)$ a structure of
enriched category $\Qcoh(G,X)$ over $\Qcoh Y$ such that
for two $G$-sheaves $\mc F$, $\mc
F'$ :
$$\Qcoh(G,X)(\mc F, \mc F')=\pi_*^G(\Qcoh X(\mc F, \mc F'))$$

The functor $\pi_*^G :\qcoh(G,X) \ra \qcoh Y$ can be lifted
 in a natural way as a  $\Qcoh Y $-functor
 $\pi_*^G :\Qcoh(G,X) \ra \Qcoh Y $.
\end{lem}

\begin{proof}
 We can first define a closed category $\Qcoh(G,X)$ whose underlying
 category is $\qcoh(G,X)$ by setting $\Qcoh(G,X)(\mc F, \mc F')
 =\Qcoh X(\mc F, \mc F')$, seen
 with its usual action, and taking as ``structural functor'' $H^0(X,\cdot)^G$.
Moreover, the functor $\pi_*^G :\qcoh(G,X) \ra \qcoh Y$ can be lifted
 in a natural way to a closed functor
 $\pi_*^G :\Qcoh(G,X) \ra \Qcoh Y $, and
 applying \cite{EilKel}, Chapter I, Proposition 6.1,
 we get the wished structure.
\end{proof}

In the sequel, we will not use the closed structure of $\Qcoh(G,X)$
described in the proof above, so we will always use this notation to
refer to the (poorer) structure of $\Qcoh Y$-category.

\begin{defi}
Let $X$ be a $G$-scheme over $k$,
$s_X:X \ra \spec k$ the structure morphism, and $\mc A$ a full
subcategory of $k[G]\,\rm \bf {mod}$.

The \emph{Auslander algebra over $Y$} associated is the ring $\mc A_X$
over $Y$ equal to the full subcategory of $\Qcoh(G,X)$ whose objects are
of the form $s_X^* V$, for all objects $V$ of $\mc A$.

In other words the objects of $\mc A_X$ are the same as those of $\mc
A$, the morphism sheaves are given by :
$\mc A_X(V,V') = \Qcoh(G,X)(s_X^* V, s_X^* V')$,
and the neutral and the composition are those induced by the ones of
$\Qcoh(G,X)$.

\end{defi}

\subsubsection{Comparison to the constant algebra}

\begin{defi}
Let $X$ be a $G$-scheme over $k$, $s_Y:Y\ra \spec k$ the structure
morphism of the quotient, and $\mc A$ a full subcategory of
$k[G]\,\rm \bf {mod}$.

The \emph{constant Auslander algebra over $Y$} is the ring
$\mc A_X^c=s_Y^*\mc A$.
\end{defi}

\begin{prop}
There is a natural morphism of rings over $Y$ : $ \mc A_X^c \ra \mc A_X$
which is an isomorphism if the action of $G$ on $X$ is trivial.
\end{prop}

\begin{proof}
To define the morphism, we can use adjunction, and the second
assertion is clear.
\end{proof}

\subsubsection{Comparison to the ''free'' algebra}

\begin{defi}
Let $X$ be a $G$-scheme over $k$,
$s_X:X \ra \spec k$ the structure morphism, and $\mc A$ a full
subcategory of $k[G]\,\rm \bf {mod}$.

The \emph{ ''free'' Auslander algebra over $Y$}
associated is the ring $\mc A_X^f$ over $Y$ equal to the full
subcategory of $\Qcoh Y$ whose objects are
of the form $\pi_*^G s_X^* V$, for all objects $V$ of $\mc A$.

In other words the objects of $\mc A_X^f$ are the same as those of $\mc
A$, the morphism sheaves are given by :
$\mc A_X^f(V,V') = \Qcoh Y(\pi_*^G s_X^* V, \pi_*^G s_X^* V')$,
and the neutral and the composition are those induced by the ones of
$\Qcoh Y$.
\end{defi}

\begin{prop}
There is a natural morphism of rings over $Y$ : $ \mc A_X \ra \mc A_X^f$
which is an isomorphism if the action of $G$ on $X$ is free.
\end{prop}

\begin{proof}
The existence of the morphism is a no more than the last assertion of
Lemma \ref{QcohGX}. When the action is free, it is a classical
result in descent theory that $\pi_*^G : \qcoh(G,X) \ra \qcoh Y$ is
an equivalence of categories, and this implies the Proposition.
\end{proof}

\subsubsection{Functoriality}

\label{FunctRinged}

Let $f : X' \ra X $ a morphism of $G$-schemes.
This defines a map between quotient schemes fitting
in a commutative diagram~:

\begin{equation}
\label{diagquot}
\xymatrix {
   X'\ar[d]_{\pi'}\ar[r]^f &      X \ar[d]^{\pi}\\
   Y' \ar[r]_{\tilde f}    &  Y
   }
\end{equation}

\begin{lem}
\label{MorG->MorA}
(i) $f : X' \ra X $ induces a morphism of ringed schemes $(\tilde {f},
  \tilde {f}^\#):(Y', \mc
  A_{X'}) \ra (Y, \mc A_{X})$.

(ii)  $adj(\tilde {f}^\#) : \tilde {f}^*  \mc A_{X} \ra \mc A_{X'}$
  is an isomorphism if the diagram
\ref{diagquot} above is fibred, i.e. if $X' = X \times_Y Y'$.

\end{lem}

\begin{proof}
(i) To construct $\tilde {f}^\#$, we start from the isomorphism
  $f^*\Qcoh X({s_X}^*V, {s_X}^*W) \ra \Qcoh X'({s_{X'}}^*V, {s_{X'}}^*W)$
  given by the fact that $f^*:\Qcoh X \ra \Qcoh X'$ is a closed
  functor. This is in fact a $G$-isomorphism and give by adjunction
  a $G$-morphism
  $\Qcoh X({s_X}^*V, {s_X}^*W) \ra f_*\Qcoh X'({s_{X'}}^*V,
  {s_{X'}}^*W)$. Applying $\pi_*^G$ and using the fact that
  $\pi_*^G f_* = \tilde{f}_*{\pi '}_*^G  $, we get the map $\mc A_X
  ({s_X}^*V, {s_X}^*W) \ra \tilde{f}_* \mc A_{X'} ({s_{X'}}^*V, {s_{X'}}^*W)$
  that we needed.

(ii) By base change, the canonical $2$-arrow $\tilde{f}^* \pi_*^G
  \Longrightarrow {\pi '}_*^G f^*$ is an isomorphism, and the result
  follows.
\end{proof}

\subsubsection{Change of group}

\label{ParaChangeGroup}
To deal with the problem of change of group, we have to enlarge
slightly our definition of the Auslander algebra to include the case
of the basic data being a \emph{functor} $F: \mc A \ra k[G]\,\mod$, not
only an inclusion. We still denote by $\mc A_X$ (instead of the better
$\mc A_F$) the corresponding ring over $Y=X/G$ whose objects are those
of $\mc A$, and whose morphisms are given by $\mc A_X(V,V') =
\Qcoh(G,X)(s_X^* FV, s_X^* FV')$.

Now let $\alpha : H \ra G$ be a group morphism. We define $(\alpha^*\mc
A)_X$ as the ring on $Z=X/H$ corresponding to the functor $\alpha^*F
: \mc A \ra k[G]\,\mod \ra k[H]\,\mod$. Let $\tilde{\alpha} : Z \ra
Y$ be the canonical morphism.
There is a natural morphism ${\tilde{\alpha}}^{\#}:\mc A_X \ra
\tilde{\alpha}_*(\alpha^*\mc A)_X$ of rings over $Y$, in other words,
we have a morphism of ringed schemes
$$(\tilde{\alpha},{\tilde{\alpha}}^{\#}): (Z,(\alpha^*\mc A)_X)\ra (Y,
\mc A_X)$$

Suppose moreover that $\alpha$ is an inclusion, and $X=G\times ^H X'$,
for a $H$-scheme $X'$, with quotient $Y'=X'/H$. We have a canonical
$H$-morphism $X' \ra X_{|H}$ and a corresponding morphism of ringed
schemes for $\alpha^*\mc A$.

\begin{lem}
\label{ChangeGroup}
With notations as above the canonical morphism of ringed schemes
$$(Y',(\alpha^*\mc A)_{X'}) \ra (Z,(\alpha^*\mc A)_X)\ra (Y,
\mc A_X)$$ is an isomorphism.
\end{lem}

\begin{proof}
This is the fact that restriction along $X' \ra X_{|H} \ra
X$ induces an equivalence $\qcoh(G,X) \simeq \qcoh(H,X')$.
\end{proof}

\subsection{$\mc A$-sheaves}

\subsubsection{Definition}

\begin{defi}
Let $X$ be a $G$-scheme over $k$, $\pi:X\ra Y$ the quotient, and
$\mc A$ a full subcategory of $k[G]\,\rm \bf {mod}$.
A $\mc A$-sheaf on $X$ is, by definition, a quasicoherent sheaf for
the ringed space $(Y, \mc A_X)$. More precisely we define $\Qcoh(\mc
A, X)$ as $\Qcoh(Y,\mc A_X)$, and $\qcoh(\mc A, X)$ as $\qcoh(Y,\mc A_X)$.
\end{defi}

\subsubsection{From $G$-sheaves to $\mc A$-sheaves}

\label{Gsh>Ash}

For any $\mc A$, we have a $\Qcoh Y$-functor :

\begin{equation*}
 \xymatrix @R=0pt{
   U_{\mc A}:\Qcoh( G,X) \ar[r]  &  \Qcoh(\mc A,X) & \\
          \mc F        \ar[r] &  \udl{\mc F} = (V \ar[r] & \Qcoh(G, X)(s_X^*V, \mc F))
 }
\end{equation*}
obtained by composing the Yoneda embedding
$\Qcoh( G,X) \ra [\Qcoh( G,X)^{op}, \Qcoh Y]$ and restriction along
$\mc A_X^{op} \ra \Qcoh( G,X)^{op}$.

\begin{lem}
\label{OneObject}
Let $\mc A=\langle k[G]  \rangle$ be the category with only the free
object $k[G]$.
Then $ U_{\langle k[G]  \rangle}$ is an equivalence of categories.
\end{lem}

\begin{proof}
Notice that $\langle k[G]  \rangle_X$ is defined as the sheaf of one object
algebras on $Y$ given by $\Qcoh(G,X)(s_X^*(k[G]),s_X^*(k[G]))$, but
this is easily seen as isomorphic to $((\pi_* \mc O_X) * G)^{op}$,
the opposite of the sheaf of twisted algebras defined by the action of
$G$ on $\pi_*\mc O_X$. So $\langle k[G]  \rangle_X^{op}$ is identified with
$(\pi_* \mc O_X) * G$, and under this isomorphism
$ U_{\langle k[G]  \rangle}$ sends
the $G$-sheaf $\mc F$ to the $(\pi_* \mc O_X) * G$-sheaf $\pi_*\mc
F$. Since $\pi$ is affine, this is an equivalence.
\end{proof}

\begin{prop}
\label{Gsheaves->Asheaves}
Let $X$ be a $G$-scheme over $k$, $\mc A$ a full
subcategory of $k[G]\,\rm \bf {mod}$ containing
$k[G]$, the free object of rank $1$.
Then $\Qcoh( G,X)$ is a reflective subcategory of $\Qcoh( \mc
A,X)$. More precisely the functor $U_{\mc A}$ admits a left adjoint $R$ such
that the counit $R U_{\mc A} \Longrightarrow 1 $  is an isomorphism.
\end{prop}

\begin{proof} The only thing to check to be able to apply Corollary
  \ref{EnrRightKan} is that the composite of the equivalence
  $U_{\langle k[G]  \rangle}$ and the right
  Kan extension $\Qcoh(\langle k[G]  \rangle,X) \ra \Qcoh(\mc A,X)$ coincides with
  $U_{\mc A}$, but this is immediate.
\end{proof}

\begin{cor}
Suppose moreover that $\mc A$ contains only projective objects. Then
$U_{\mc A}$ is an equivalence.
\end{cor}

\begin{proof}
Lemma \ref{EnrRightKan2} allows to reduce to the case where $\mc
A=\langle k[G]  \rangle$, which was the object of Lemma \ref{OneObject}.
\end{proof}

\subsubsection{Change of ring}

\begin{prop}
\label{EnrRightAdj}
\ \\
Let $X$ be a $G$-scheme over $k$, $\mc A$ a full
subcategory of $k[G]\,\rm \bf {mod}$,
and $\mc A'$ be a full subcategory of $\mc A$.

(i) The restriction functor $R:\Qcoh(\mc A,X)
\ra \Qcoh(\mc A',X) $ admits as a right adjoint the functor $K$
defined on objects by :

\begin{equation*}
 \xymatrix @R=0pt {
  K: \Qcoh(\mc A',X) \ar[r]  &  \Qcoh(\mc A,X) & \\
          \mc F        \ar[r] &  (V \ar[r] & \Qcoh(\mc A',X)
(\mc A_X(\cdot,V)_{|\mc A'_X },\mc F))}
\end{equation*}
Moreover $RK \simeq 1$ (equivalently, $K$ is fully faithful).

(ii) In particular, if $\mc A$ is projectively complete,
and every object of $\mc A$ is a direct summand of a
finite direct sum of objects of $\mc A'$, this adjunction is an
equivalence $\Qcoh(\mc A,X) \simeq  \Qcoh(\mc A',X)$.
\end{prop}

\begin{proof}
(i) follows from Corollary \ref{EnrRightKan} and
Proposition \ref{QcohCompl}, and (ii) from
Corollary \ref{EnrRightKan2}.
\end{proof}

\subsubsection{Action of $\Qcoh Y$}

\label{action}

Starting from the external ${\mc H om}$ :

\xymatrix@R=0pt{
  &&(\Qcoh Y)^{op} \otimes \Qcoh(\mc A, X) \ar[r] & \Qcoh(\mc A, X)& \\
  && \mc G \otimes  \mc F  \ar[r] & {\mc H om}(\mc G, \mc F) =(V
  \ar[r] & {\mc H om}(\mc G, \mc F(V))) }

we get as a left adjoint an action of $\Qcoh Y$ on $\Qcoh(\mc A, X)$
:

\xymatrix@R=0pt{
  &&\Qcoh Y \otimes \Qcoh(\mc A, X) \ar[r] & \Qcoh(\mc A, X) \\
  && \mc G \otimes  \mc F  \ar[r] &  \mc G \otimes  \mc F}

\subsubsection{Functoriality}

\begin{defi}
Let $f : X' \ra X $ be a morphism of $G$-schemes, and
$(\tilde {f}, \tilde {f}^\#)$ the
morphism of ringed schemes associated by Lemma \ref{MorG->MorA}.
The pull-back $f^{\triangle}$ (resp. the push-forward $f_{\triangle}$)
given by Definition \ref{FunctRingedSchemes} will be denoted by
$f^{\mc A}$ (resp. $f_{\mc A}$).
\end{defi}

In particular, for each point $Q$ of $\mc A$, we get a stalk functor :

$$\xymatrix@R=0pt{
  &&& \Qcoh(\mc A, X) \ar[r] & \Qcoh(\mc A, X\times_Y {\rm spec}\, \mc O_{Y,Q}) \\
  &&& \mc F \ar[r] & \mc F_Q }$$

Indeed, with the notations of \S \ref{StalksRinged}, Lemma
\ref{MorG->MorA} (ii) implies that $(\mc A_X)_Q \simeq \mc A_{X\times_Y {\rm spec}\, \mc O_{Y,Q}}$.

\subsubsection{Adjunction}

\begin{prop}
\label{AdjunctionAsheaves}
Let  $f : X' \ra X $ be a morphism of $G$-schemes.
The couple $(f^{\mc A},f_\mc A)$ is part of a $\Qcoh Y$-adjunction
between $\Qcoh(\mc A, X)$ and $\tilde{f}_*\Qcoh(\mc A, X')$.
\end{prop}

\begin{proof}
This follows from Proposition \ref{AdjunctionRingedSchemes}.
\end{proof}

\subsubsection{Representable sheaves}

Let $X$ be a $G$-scheme, and $\mc A$ be a full subcategory of
$k[G] \mod$. Then the category $\Qcoh(\mc A,X)=[\mc
A_X^{op},\Qcoh Y]$ contains for each $V$ the corresponding
representable functor, which we denote by  $(\mc A_X)_V$, and call a
\emph{representable sheaf}. There is an obvious local notion of
\emph{locally representable sheaf}. Moreover, if $f:X'\ra X$ is any
$G$-morphism, then one checks that $f^{\mc A}(\mc A_X)_V \simeq (\mc
A_{X'})_V$, hence both notions are preserved by arbitrary pullback.

\subsubsection{Change of group}

We keep the notations of \S \ref{ParaChangeGroup} : $X$ is a
$G$-scheme, $\alpha : H \ra G$ a group morphism, and $F: \mc A \ra
k[G] \mod$ a functor.
The associated morphism of ringed schemes
$$(\tilde{\alpha},{\tilde{\alpha}}^{\#}): (Z,(\alpha^*\mc A)_X)\ra (Y,
\mc A_X)$$ and Definition \ref{FunctRingedSchemes} provides
a restriction functor
$\alpha^{\mc A} : \Qcoh(\mc A, X) \ra
 \tilde{\alpha}_* \Qcoh(\alpha^*\mc A, X_{|H})$
and an induction functor
$\alpha_{\mc A} : \tilde{\alpha}_*\Qcoh(\alpha^*\mc A, X_{|H})
\ra \Qcoh(\mc A, X)$.
Again, Proposition \ref{AdjunctionRingedSchemes} implies that
 $(\alpha^{\mc A}, \alpha_{\mc A})$ is part of a natural adjunction
 between $\Qcoh(\mc A, X)$ and $\tilde{\alpha}_* \Qcoh(\alpha^*\mc A, X_{|H})$.

If moreover that $\alpha$ is an inclusion, and $X=G\times ^H X'$,
for a $H$-scheme $X'$, with quotient $Y'=X'/H$, Lemma \ref{ChangeGroup}
implies that there is a canonical equivalence
$$\Qcoh(\mc A, X) \simeq \Qcoh(\alpha^*\mc A, X')$$

\subsubsection{Internals homs and tensor product}

\label{Amonoidal}

For the moment being, we have just considered $\Qcoh(\mc A, X)$ as an
enriched category over $\Qcoh Y$. But the work of B.Day (see
\cite{Day}) implies that if $\mc A$ is a submonoidal category of
$k[G] \mod$ (endowed with the tensor product over $k$), then
$\Qcoh(\mc A, X)$ carries the structure of a monoidal closed symmetric
category, for which the functor of evaluation at the unit
$\cdot(k) : \Qcoh(\mc A, X) \ra \Qcoh Y$ is closed. Since by pushing
$\Qcoh(\mc A, X)$, seen as enriched over itself, along this functor,
we recover $\Qcoh(\mc A, X)$ seen as enriched category over
$\Qcoh Y$, we keep the same notation.

The starting fact is the following : suppose $\mc A$ is a
full submonoidal category of $k[G] \mod$, and $X$ is a $G$-scheme.
Then the Auslander algebra $\mc A_X^{op}$ has a natural structure of a
monoidal symmetric category \emph{over $\Qcoh Y$}. So \cite{Day} \S
3,\S4 shows that there is a canonical structure of monoidal closed symmetric
category on $\Qcoh(\mc A,X)$, whose unit object is

$$ \udl{\mc O_X} : V \ra \pi_*^G(s_X^*V^{\vee})$$
whose internal homs are given by :

$$\Qcoh(\mc A,X)(\mc F, \mc G)(V)=\int_{W}\Qcoh Y(\mc F(W),\mc
G(W\otimes_k V))$$
and whose tensor product is given by a convolution formula :

$$ \mc F \otimes \mc G (V) = \int^W  \mc F(W) \otimes_{\mc O_Y} \mc
G(V\otimes_k W^{\vee})$$

Tensor product with a representable sheaf can be made more explicit :

\begin{lem}
\label{tensrepr}
$$ (\mc F \otimes (\mc A_X)_V)(W) \simeq \mc F ( W\otimes_k V^{\vee})$$
\end{lem}

\begin{proof}
Since $\mc A_X$ is dense (or adequate) in $\Qcoh(\mc A,X)$, it
suffices to check this on $\mc F =(\mc A_X)_{V'}$. But then it boils
down to the fact that the tensor product on $\Qcoh(\mc A,X)$ extends
the tensor product on $\mc A_X$.
\end{proof}

We deduce a projection formula in this context :

\begin{prop}
Let $f:X'\ra X$ be any $G$-morphism, $\mc F$ a locally representable
$\mc A$-sheaf on $X$, $\mc G$ a quasicoherent $\mc A$-sheaf on $X'$.
Then the natural morphism
$$ \mc F \otimes f_{\mc
  A}\mc G \ra f_{\mc A}(f^{\mc A}\mc F \otimes \mc G)$$
is an isomorphism.
\end{prop}

\begin{proof}
This is a local problem, so we can suppose $\mc F$
representable. But now we can use Lemma \ref{tensrepr} to conclude.
\end{proof}

\subsubsection{Cohomology}
\label{Cohomology}
\begin{prop}
Let $X$ be a $G$-scheme over $k$, and $\mc A$ a full subcategory
of $k[G]\mod$. The category $\qcoh(\mc A,X)$ has enough injective
objects.
\end{prop}

\begin{proof}
Using Propositions \ref{ExactStalks} and \ref{AdjunctionAsheaves}, we see
that the classical proof (see \cite{Tohoku}) applies without change.
\end{proof}

Proposition \ref{AdjunctionAsheaves} also shows that the functor $f_{\mc
A}$ is left exact, hence the following definition.

\begin{defi}
Let $f: X' \ra X$ be a morphism of $G$-schemes. We denote by
$R^if_\mc A : \qcoh(\mc A, X') \ra \qcoh(\mc A, X)$ the $i$-th
derived functor of $f_{\mc A}$.
\end{defi}

In particular, when $X={\rm  spec}\, k$, we denote the derived functors of
the global sections functor by $H^i(X',\cdot)$.

In view of Corollary \ref{Gsheaves->Asheaves}, it is natural to
compare the usual cohomology of a $G$-sheaf to the one of the corresponding
$\mc A$-sheaf. We give three comparison results.

\begin{prop}
\label{precomp}

Let $X$ be a $G$-scheme over $k$, $\mc A$ a full
subcategory of $k[G]\,\rm \bf {mod}$ containing
$k[G]$, the free object of rank $1$.
Suppose given a $G$-sheaf $\mc F$ on $X$. There is a canonical
$G$-isomorphism : $H^i(X,\udl {\mc F})(k[G]) \simeq H^i(X, {\mc F})$
\end{prop}

\begin{proof}
This an immediate consequence of the exactness of the evaluation
$\cdot (k[G])$ and of the isomorphism $\udl {\mc F} (k[G]) \simeq
\pi_* \mc F$.
\end{proof}

\begin{prop}
\label{comp}
Let $X$ be a $G$-scheme over $k$, $\mc A$ a full
subcategory of $k[G]\,\rm \bf {mod}$ containing
$k[G]$, the free object of rank $1$.
Suppose that the action of $G$ on $X$ is tame.
Then for each $G$-sheaf on $X$ we have a spectral sequence :

$${\rm Ext}^p_{k[G]}(\cdot,H^q(X,\mc F))\Rightarrow H^{p+q}(X,\udl {\mc F})$$
\end{prop}

\begin{proof}
The tameness of the action says that the functor $\pi_*^G$ is exact,
and so is the functor $\mc F \ra \udl {\mc F}$. Hence the result is a
direct consequence of a Theorem of Grothendieck describing the derived
functors of a composite functor.
\end{proof}

In particular, both cohomology coincide for a reductive action, i.e.
when the characteristic of $k$ does not divide the order of $G$. But
note that they can differ even for a free action, as soon as $p={\rm
  car}k | \#G$.

Remember that the equivariant cohomology functors $H^i(X,G,\cdot)$ are
defined as the derived functors of $H^0(X, \cdot)^G$ (see \cite{Tohoku}).

\begin{prop}
Let $X$ be a $G$-scheme over $k$, $\mc A$ a full
subcategory of $k[G]\,\rm \bf {mod}$ containing
$k[G]$, the free object of rank $1$, and $k$, the trivial representation.
Then for each $G$-sheaf on $X$ we have a spectral sequence :
$$H^p(X, R^q U_{\mc A} \mc F)(k)\Rightarrow H^{p+q}(X,G,\mc F)$$
\end{prop}

\begin{proof}
This is a consequence of the isomorphism
$H^0(X,\udl {\mc F})(k) \simeq H^0(X, \mc F)^G$.
\end{proof}

\subsubsection{Canonical filtration}

\label{CanFiltr}

Carrying on with the idea of S.Nakajima in \cite{Naka},
we will see that, if we impose some
more structure on $\mc A$, then the kernel preserving $\mc A$-sheaves
(and in particular those coming from $G$-sheaves) have a useful natural
filtration of algebraic nature.

Since $k[G]$ is a semilocal ring, for each object $V$ of $k[G] \mod$,
we have ${\rm rad} V =({\rm rad} \; k[G]) V$. In particular there is a
functor ${\rm rad} : k[G] \mod \ra k[G] \mod$. Moreover we have a
natural duality functor $D:   k[G] \mod \ra (k[G] \mod)^{op}$ sending
$V$ to its dual $V^{\vee}$.

The data we need is the following : $\mc A$ is as usual a full subcategory of
$k[G] \mod$, that we suppose stable under ${\rm rad}$ and $D$. We need
also a compatibility between these two operations, in the sense that
we suppose given a natural isomorphism $\alpha$~:

\xymatrix@R=12pt{
  &&&& \mc A \ar[rr]^{D} \ar[dd]_{rad} &&\mc A^{op}
   \ar[dd]^{rad^{op}}\\
  &&&&&&\\
  &&&& \mc A \ar[rr]_{D} \ar@{<=}[uurr]^{\alpha}&&\mc A^{op} }

Using $\alpha$, the natural transformation $rad \Longrightarrow 1$,
and the canonical isomorphism $\beta : D \; D^{op} \Longrightarrow1$
we get finally a natural transformation ${\rm rad}^{op}
\Longrightarrow 1$, as sketched in the following diagram :

\xymatrix{
&&& \mc A^{op} \ar[r]^{D^{op}} \ar@/^3pc/[rrr]|*{}="BAS2"_{DD^{op}}
  \ar@/^4pc/[rrr]|*{}="HAUT2"^1 \ar @{=>} "BAS2";"HAUT2"
  \ar@/_3pc/[rrr]|*{}="HAUT3"^{DradD^{op}}
  \ar@/_4pc/[rrr]|*{}="BAS3"_{rad^{op}} \ar @{=>} "BAS3";"HAUT3"
  & \mc A \ar@/^/[r]|*{}="HAUT"^1  \ar@/_/[r]|*{}="BAS"_{rad}
\ar @{=>} "BAS";"HAUT" & \mc A \ar[r]^D & \mc A^{op} }
\noindent where the bottom natural transformation is $\alpha D^{op} \circ
rad^{op} \beta ^{-1}$.

This is a priori a natural transformation between additive functors
from $\mc A^{op}$ to  $\mc A^{op}$, but if $X$ is a $G$-scheme, this
extends to a natural transformation of the corresponding functors from
$\mc A_X^{op}$ to  $\mc A_X^{op}$, which we note the same way.

Applying now the $2$-functor : $(\cdot,\Qcoh Y): (\Qcoh
Y-CAT)^{op} \ra {\mc C}at$ we get a natural transformation $({\rm
  rad}^{op},\Qcoh Y) \Longrightarrow 1$ between endofunctors of
$\qcoh(\mc A,X)$. We define $R$ as  $({\rm rad}^{op},\Qcoh Y) $.

Let now $\mc F$ be a $\mc A$-sheaf on $X$ such that $\mc F$, as a
functor, preserves monomorphisms (this is in particular the case if $\mc
F=\udl{\mc G}$ for a $G$-sheaf $\mc G$). Then the natural morphism
$R \mc F \ra \mc F$ is itself a monomorphism, and moreover $R\mc F$ preserves
monomorphisms, so that $\mc F$ has a canonical filtration. We sum up
the construction in the following definition :

\begin{defi}
Let $\mc A$ be a full subcategory of
$k[G] \mod$, stable under radical $rad$ and duality $D$.

(i)
We will say that \emph{radical and duality commute} to say that we fix
an isomorphism ${\rm rad}^{op} D \Longrightarrow
D rad $.

(ii) Suppose that radical and duality commute, and let
moreover $X$ be a $G$-scheme, and $\mc F$ a $\mc A$-sheaf
on $X$, whose underlying functor preserves monomorphisms.
In this circumstances, we will denote the induced filtration on $\mc
F$ by

$$\cdots \subset R^i \mc F \subset \cdots  \subset R^2\mc F \subset R^1\mc F
\subset \mc F$$
and the associated graded $\mc A$-sheaf by
$${\rm gr}_i \mc F =  R^i \mc F /R^{i+1} \mc F$$

\end{defi}

The existence of an isomorphism ${\rm rad}^{op} D \Longrightarrow
D rad $ is not the general rule. Indeed, we have  a natural
isomorphism :

$$ rad(V^\vee) \simeq (\frac{V}{soc V})^{\vee} $$
where $soc V$ is the sum of all simple submodules of $V$,
but in general $V/soc V$ is \emph {not} isomorphic to $rad V$.

We finish by a positive example, to be used in the sequel : let $G$ be
a cyclic of order $n$, with fixed generator $\sigma$. We choose as
ring with several object $\mc A = k[G] \, \mod$, the whole category
of all $k[G]$-modules of finite type.

The representation theory of $G$ is easily described, as follows.
Write $n=p^v a$, where $a$ is prime to $p$. For each integer $0
\leq j \leq p^v$ define $V_j = k[G]/(\sigma -1)^j$. Then the set
of modules $\{ \psi \otimes V_j \; / \; \psi \in \widehat G, \; 1
\leq j \leq p^v \}$ is a skeleton of the indecomposables in $\mc
A$.

To construct a natural isomorphism $V/soc V \simeq rad V$, we can
restrict to the indecomposable skeleton we have just fixed. Now the
maps
$$ \psi \otimes (\sigma -1) : \psi \otimes V_j  \ra  \psi \otimes V_j
$$
give a natural transformation, $1 \Longrightarrow 1 $, which
canonically factorizes in a natural isomorphism $V/soc V \simeq rad V$.


\section{{Modular $K$-theory}}

\label{Modular}

\subsection{Definition}

\begin{defi}
Let $X$ be a $G$-scheme over $k$ and $\mc A$ a full
subcategory of $k[G]\,\rm \bf {mod}$.
A quasicoherent
$\mc A$-sheaf $\mc F$ on $X$ is said \emph{coherent} if for each
$G$-invariant
open affine $U = {\rm spec}\, R$ of $X$, the restriction $\mc F_{|U}$ is
of finite type in $\Qcoh(\mc A, U) \simeq [\mc A_X(U)^{op}, R^G \,
\Mod]$ (i.e. if, seen in $[\mc A_X(U)^{op}, R^G \,
\Mod]$, it is a quotient of a finite sum of representable objects).
We denote by $\coh(\mc A, X)$ the whole subcategory of $\qcoh(\mc A,
X)$ whose objects are the coherent $\mc A$-sheaves.
\end{defi}

This notion is Morita invariant, at least in the following sense :

\begin{prop}
\label{CohMor}
Let $X$ be a $G$-scheme over $k$, $\mc A$ a full
subcategory of $k[G]\,\rm \bf {mod}$,
and $\mc A'$ be a full subcategory of $\mc A$.
Suppose that $\mc A$ is projectively complete,
and that every object of $\mc A$ is a direct summand of a
finite direct sum of objects of $\mc A'$. Then restriction along
$\mc A' \ra \mc A$ induces an
equivalence $\coh(\mc A,X) \simeq  \coh(\mc A',X)$.
\end{prop}

\begin{proof}
We know from Proposition \ref{EnrRightAdj} that restriction along
$\mc A' \ra \mc A$ induces an equivalence $\qcoh(\mc A,X) \simeq
\qcoh(\mc A',X)$. Because of the hypothesis, this restriction sends
$\coh(\mc A,X)$ to $\coh(\mc A',X)$. Moreover the left adjoint
$\otimes_{\mc A'_X}\mc A_X$ is an inverse equivalence, and since it is
right exact and preserves representables, it sends $\coh(\mc A',X)$ to
$\coh(\mc A,X)$.
\end{proof}

The $\Qcoh Y$ enriched functor $\Coh Y \ra \Qcoh Y$ allows to
identify $\Qcoh Y-CAT(\mc A_X^{op}, \Coh Y)$
with a subcategory of $\qcoh(\mc A,X)$.
Since $\mc A_X$ is in fact enriched in $\Coh Y$, $\coh(\mc A,
X)$ is a subcategory of $\Qcoh Y-CAT(\mc A_X^{op}, \Coh Y)$. Moreover~:

\begin{lem}
  \label{CohCoh2}
Suppose $\mc A$ admits a
 finite set of additive generators.
Then $\coh(\mc A, X)=\Qcoh Y-CAT(\mc A_X^{op}, \Coh Y)$.
\end{lem}

\begin{proof}
This is a local question, hence we can conclude by applying Lemma
\ref{CohCoh}, which is of course valid with the field $k$ replaced
by any commutative noetherian ring.
\end{proof}

\begin{lem}
Suppose $\mc A$ admits a
 finite set of additive generators.
Then $\coh(\mc A, X)$ is an abelian category.
\end{lem}

\begin{proof}

This is clear from Proposition \ref{CohMor},
which allows to reduce to the case when $\mc A$ has only one object.
This follows also from Lemma \ref{CohCoh2}, because we can follow word for
word the proof of Proposition \ref{abelian}.
\end{proof}

\begin{defi}
Let $X$ be a $G$-scheme over $k$ and $\mc A$ a full
subcategory of $k[G]\,\rm \bf {mod}$  admitting
a finite set of additive generators.
We denote by $K_i(\mc A, X)$ the Quillen $i$-th group of the abelian
category $\coh(\mc A, X)$.
\end{defi}

\subsection{Functoriality}

\subsubsection{Pullback}

Let $f: X'\ra X$ be a morphism of $G$-schemes over $k$ such that the
morphism $\tilde{f} : Y' \ra Y$ between quotient schemes is flat, and
$X'=X\times_Y Y'$.
Then the functor $f^{\mc A}:\coh(\mc A, X) \ra \coh(\mc A, X')$ is
exact, hence induces a map in $K$-theory.

\subsubsection{Pushforward}

\begin{lem}
  \label{finite}
Let $f: X'\ra X$ be a morphism of $G$-schemes over $k$ and $\mc A$ a full
subcategory of $k[G]\,\rm \bf {mod}$, admitting a
finite set of additive generators. Suppose that the morphism $\tilde f
: Y' \ra Y$ between quotient schemes is proper, then :

(i) For each coherent $\mc A$-sheaf $\mc F$ on $X'$, and each nonnegative
integer $i$, the $\mc A$- sheaf $R^if_{\mc A}\mc F$ is coherent.

(ii) There exists an integer $n$, such for any integer $i>n$, and any
coherent $\mc A$-sheaf $\mc F$ on $X'$, we have $R^if_{\mc A}\mc F=0$.
\end{lem}

\begin{proof}
Because cohomology commutes with the projection $p_V$ (i.e. for any $V$
in ${\rm obj} \mc A$, $R^if_{\mc A}\mc F(V)=R^if_{\langle V \rangle}(\mc
F(V))$),
Lemma \ref{CohCoh2} allows to reduce to the case where $\mc A$ has
only one object (one can also use Proposition \ref{CohMor} to reduce
to this case). But since $R^if_{\langle V \rangle}(\mc
F(V))$, seen in $\qcoh Y$, is nothing else that $R^i\tilde{f}(\mc
F(V))$, the Lemma results from \cite{EGA3}, 3.2.1, 1.4.12.
\end{proof}

Now, given such a $f: X'\ra X$, we can follow the argument given
in \cite{Quil} \S 7, 2.7, to define a map

$$f_{\mc A}: K_i(\mc A, X') \ra K_i(\mc A, X)$$

in the following two cases :

(i) $\tilde f$ is finite,

(ii) $Y'$ admits an ample line bundle (then we have to use the action
of $\qcoh Y$ on $\qcoh(\mc A,X)$ defined in section \ref{action}).

Given $f: X'\ra X$ and $g: X''\ra X'$, both checking the condition of
Lemma \ref{finite}, and one of the conditions above, then the formula
$(fg)_{\mc A}=f_{\mc A}g_{\mc A}$ holds.

\subsection{Localization}

\begin{thm}
\label{Loc}
  Let $i: X'\ra X$ be a morphism of $G$-schemes over $k$, and $\mc A$ a full
subcategory of $k[G]\,\rm \bf {mod}$, admitting a
finite set of additive generators.

Suppose that the morphism $\tilde i : Y' \ra Y$
between quotient schemes is a closed immersion, and that $i^\# :
\mc A_X \ra \tilde i_*\mc A_{X'}$ is an epimorphism.

Denote by $U$ the pullback by $\pi:X \ra Y$ of the complement of $Y'$
in $Y$, and by $j:U\ra X$ the canonical inclusion.

Then there is a long exact sequence :

\xymatrix @R=40pt{
 &&& \cdots & \cdots & \cdots
 \ar `dr_l[ll]  `^dr[lld]  [lld] \\
&&& K_i(\mc A, X') \ar[r]^{i_{\mc A}} & K_i(\mc A,X) \ar[r]^{j^{\mc A}} & K_i(\mc A,U)
 \ar `dr_l[ll]  `^dr[lld]  [lld] \\
 &&& \cdots & \cdots & \cdots
 \ar `dr_l[ll]  `^dr[lld]  [lld] \\
 &&& K_1(\mc A, X') \ar[r]^{i_{\mc A}} & K_1(\mc A,X) \ar[r]^{j^{\mc A}} & K_1(\mc A,U)
 \ar `dr_l[ll]  `^dr[lld]  [lld] \\
 &&& K_0(\mc A, X') \ar[r]^{i_{\mc A}} & K_0(\mc A,X) \ar[r]^{j^{\mc A}} & K_0(\mc A,U) \ar[r] &0 }

\end{thm}

\begin{proof}
The hypothesis on $\mc A$ and Proposition \ref{CohMor} allows to
reduce to the case where $\mc A$ has only one object, what we will do
from now on.

The idea is of course to apply \cite{Quil} \S 5, Theorem 5, but to do
so we need the two following facts.

Denote by $\coh(\mc A, X)_{Y'}$ the full subcategory of $\coh(\mc
A, X)$ consisting of sheaves with support with $Y'$. Being the kernel
of the restriction functor $j^{\mc A}$, this is a Serre subcategory, and the first
step consist of showing that $j^{\mc A}$ induces an equivalence

\begin{equation}
\label{AnRes}
\coh(\mc A, X)/ \coh(\mc A, X)_{Y'} \simeq \coh(\mc A, U)
\end{equation}

This is the object of section \ref{res}. Note that the notion of
quotient category
used here (quotient as an example of localization) has nothing to see
with the notion of quotient used in part \ref{Repr} (quotient by a
two-sided ideal).

Then we will show than we can apply the hypothesis of the d\'evissage
Theorem (\cite{Quil} \S 5, Theorem 4) to the functor
$i_{\mc A} : \coh(\mc A, X') \ra \coh(\mc A, X)_{Y'}$ : this is the aim of
section \ref{dev}.

We sum up the notations we use in the proof in the following diagram

\begin{equation}
\label{NotProof}
\xymatrix {
   X'\ar[d]_{\pi'}\ar[r]^i &      X \ar[d]^{\pi} & U \ar[d]^{\pi_{|U}}
     \ar[l]_j \\
   Y' \ar[r]_{\tilde i}    &  Y  & O \ar[l]^{\tilde j}}
 \end{equation}
\end{proof}

\subsection{Restriction of coherent $\mc A$-sheaves to an open}

\label{res}

Because we lack of a complete reference we recall the classical arguments.

\begin{prop}
Let $F: \mc B \ra \mc C$ be an exact functor between abelian categories
such that

(i) For any object $N$ of $\mc C$, there is an object $M$ of $\mc B$
and an isomorphism $FM \simeq N$,

(ii) For any objects $M$, $M'$ of $\mc B$, and any map $x : FM \ra FM'$ in
$\mc C$, there exists a diagram in $\mc B$

\xymatrix{
  &&&&    & M''\ar[dl]_u \ar[dr]^v & \\
  &&&&   M & & M'}

\noindent such that the diagram

\xymatrix{
  &&&&    & FM''\ar[dl]_{Fu} \ar[dr]^{Fv} & \\
  &&&&   FM\ar[rr]^x & & FM'}

\noindent commutes in $\mc C$, and such that $Fu$ is an isomorphism.

Then the canonical functor $\mc B / {\rm Ker} F \ra \mc C$ is an
equivalence of categories.

\end{prop}

\begin{proof}
(i) shows that $\mc B / {\rm Ker} F \ra \mc C$ is essentially full,
  and (ii) that it is fully faithful.
\end{proof}

Hence the equivalence \ref{AnRes} will follow from the three following
Lemmas.

\begin{lem}
 \label{key}
  Let $\mc F$ be a coherent $\mc  A$-sheaf on $X$, and $\beta : \mc G \ra
  \mc F_{|U}$ a monomorphism in $\coh(\mc A,U)$. There exists a
  monomorphism $\alpha : \mc G' \ra \mc F$ in $\coh(\mc A,X)$ such
  that $\alpha_{|U}=\beta$ (as subobjects of $\mc F_{|U}$).
\end{lem}

\begin{proof}
  As in \cite{BorSer}, Proposition 1, let $\mc G'$ be the sheaf on $Y$
  associated to the presheaf $O' \ra \{ s \in \mc F(O') / \exists
  t \in \mc G(O \cap O') / \beta(t)=s_{|O \cap O'} \}$, and
  $\alpha :\mc G' \ra \mc F$ be the canonical map. In \cite{BorSer} is
  shown that $\alpha$ is a map in $\coh Y$ such that
  $\alpha_{|U}=\beta$ as subobjects of $\mc F_{|U}$ in $\coh O$. But from its
  definition, one sees at once that $\mc G'$ is stable under the
  action of $\mc A$, hence has a unique structure of $\mc A$-sheaf
  such that $\alpha$ is an arrow in $\coh(\mc A,X)$.
\end{proof}

\begin{lem}
  Let $\mc G$ be a coherent $\mc A$-sheaf on $U$. There exists a coherent
  $\mc A$-sheaf $\mc F$ on $X$ such that $\mc F_{|U}=\mc G$.
\end{lem}

\begin{proof}
 Follows Lemma \ref{key} as in the proof of \cite{BorSer}, Proposition 2.
\end{proof}

\begin{lem}
  Let $\mc F$ and $\mc F'$ be two $\mc A$-sheaves on $X$, and $\gamma
  : \mc F_{|U} \ra \mc F'_{|U}$ a morphism in $\coh(\mc A,U)$. There
  exists a diagram in $\coh(\mc A,X)$

  \xymatrix{
  &&&&    & \mc H\ar[dl]_{\alpha} \ar[dr]^{\beta} & \\
  &&&&   \mc F & & \mc F'}

\noindent such that the diagram

\xymatrix{
  &&&&    & \mc H_{|U}\ar[dl]_{\alpha_{|U}} \ar[dr]^{\beta_{|U}} & \\
  &&&& \mc F_{|U}  \ar[rr]^{\gamma} & & \mc F'_{|U}}

\noindent commutes in $\coh(\mc A,U)$, and such that $\alpha_{|U}$ is an isomorphism.

\end{lem}

\begin{proof}
  Apply Lemma \ref{key} to the graph of $\gamma$.
\end{proof}

\subsection{D\'evissage}
\label{dev}

We have to show that each object $\mc F$ in $\coh(\mc A, X)$ with support in
$Y'$ has a finite filtration whose quotients are objects in the image
of $i_{\mc A}: \coh(\mc A, X') \ra \coh(\mc A, X)$.
Let $\mc J_{X'}= {\rm ker} (i^\# : \mc A_X \ra \tilde i_*\mc
A_{X'})$.
Since we have an exact sequence $0 \ra \mc F \mc J_{X'} \ra \mc F \ra
i_{\mc A}i^{\mc A}\mc F\ra 0 $, it is enough to show that for $n$
large enough $\mc F  \mc J_{X'}^n=0$.

Define $\mc I_{Y'}= {\rm ker} ( \mc O_Y \ra \tilde i_*\mc O_{Y'})$,
the ideal sheaf of $Y'$. Since the support of $\mc F$ is included in
$Y'$, the Nullstellensatz ensures that for $n$ large enough we have
$\mc I_{Y'}^n \subset {\rm ann}\, \mc F$, hence the following Lemma will
be enough to conclude.

\begin{lem}
For $n$ large enough $\mc J_{X'}^n \subset \mc I_{Y'} \mc A_X$.
\end{lem}

\begin{proof}
 This proof was suggested to me by A.Vistoli.

First remember that we can suppose that $\mc A$ has only one object,
and denote by $I$ this $k[G]$-module. Then $G$ acts by ring
homomorphisms on ${\rm End}_k\, I$, and if $s_X : X \ra {\rm Spec}\,
k$ denotes the structure morphism, we have that $\mc A_X= \pi_*^G
s_X^* {\rm End}_k\, I$.

Define $\mc I_{X'}= {\rm ker} ( \mc O_X \ra i_*\mc O_{X'})$.
The left exactness of $\pi_*^G$ implies that $\mc J_{X'}
=  \pi_*^G (\mc I_{X'} s_X^* {\rm End}_k\, I)$. Hence for all $n$ :

$$\mc J_{X'}^n =(\pi_*^G (\mc I_{X'} s_X^* {\rm End}_k\, I))^n
\subset \pi_*^G ((\mc I_{X'} s_X^* {\rm End}_k\, I)^n)
= \pi_*^G (\mc I_{X'}^n s_X^* {\rm End}_k\, I)$$

Fix an integer $r$. Applying the Nullstellensatz again, we have that
for large $n$ : $\mc I_{X'}^n \subset  (\pi^* \mc I_{Y'}) ^r$, hence
for large $n$ :

$$\mc J_{X'}^n \subset  \pi_*^G ((\pi^* \mc I_{Y'}) ^r
s_X^* {\rm End}_k\, I)=
\mc (I_{Y'} ^r \pi_* (s_X^* {\rm End}_k\, I))^G
=\mc I_{Y'} ^r \pi_* (s_X^* {\rm End}_k\, I) \cap \mc A_X
$$
Since $Y$ is noetherian, we can apply the Artin-Rees Lemma to conclude
that the filtration $(\mc I_{Y'} ^r \pi_* (s_X^* {\rm End}_k\, I) \cap
\mc A_X)_{r\geq 0}$ is $\mc I_{Y'}$-stable. In particular for large $r$ we
have $\mc I_{Y'} ^r \pi_* (s_X^* {\rm End}_k\, I) \cap \mc A_X \subset
\mc I_{Y'}\mc A_X$, and the Lemma is shown.

\end{proof}

\subsection{A criterion of surjectivity}

\begin{prop}

\label{CritSurj}
   Let $i: X'\ra X$ be a morphism of $G$-schemes over $k$, such that

   (i) There is a normal subgroup $H$ of $G$, such that $H$ acts
   trivially on $X'$, and $G/H$ acts freely on $X'$,

   (ii) The morphism $\tilde i : Y' \ra Y$
   between quotient schemes is a closed immersion.

   Then the canonical morphism
   $i^\# :\mc A_X \ra \tilde i_*\mc A_{X'}$ is an epimorphism.
\end{prop}

\begin{proof}
For the moment being, we do not use the hypothesis (i), and define
$\nu :X \ra Z=X/H$ and $\mu : Z \ra Y=Z/P$ as the quotient morphisms,
and similarly for $X'$, so that we get a $P$-morphism $\hat f :Z'\ra Z$
fitting in the following commutative diagram :

\begin{equation*}
\xymatrix {
   X'\ar[d]_{\nu'}\ar[r]^f &      X \ar[d]^{\nu}\\
   Z'\ar[d]_{\mu'}\ar[r]^{\hat f} &      Z \ar[d]^{\mu}\\
   Y' \ar[r]_{\tilde f}    &  Y
   }
\end{equation*}

Define a new ring $\mc A_X^{c,H}$ on $Y$ by setting

$$\mc A_X^{c,H}(V,W)= \mu_*^P(\Qcoh Z (s_Z^*V, s_Z^*W)^H)$$

where $s_Z : Z \ra {\rm spec}\, k$ is the structure morphism.

This definition is functorial in $X$.

Moreover, since $\pi_*^G =\mu_*^P \nu_*^H $, there is a canonical morphism
$\mc A_X^{c,H} \ra \mc A_X$, also functorial in $X$, so that we
get a commutative diagram :

\begin{equation*}
\xymatrix{
  \mc A_X \ar[r] &  \tilde f _*\mc A_{X'}\\
  \mc A_X^{c,H} \ar[r]\ar[u] &  \tilde f _*\mc A_{X'}^{c,H}\ar[u] }
\end{equation*}

The first part of hypothesis {\it (i)} means that $\nu'=1$, and this implies
that $\mc A_{X'}^{c,H} \simeq \mc A_{X'}$.

Since the question is local, we can use the second part of hypothesis
{\it (i)} to drop the ramification locus of $\mu$ and thus reduce to the
case where $P$ acts freely on $Z$. But then descent theory implies
that  $  \tilde f^* \mc A_X^{c,H} \simeq \mc A_{X'}^{c,H}$.

So the Proposition now follows from hypothesis {\it (ii)}.

\end{proof}

\subsection{The case of a free action}

\begin{prop}
  \label{free=>equivalence}
Let $X$ be a $G$-scheme over $k$ and $\mc A$ a full
subcategory of $k[G]\mod$  admitting
a finite set of additive generators, and containing $k[G]$, the free
object of rank $1$.
If the action of $G$ on $X$ is free, then the functor
$U_{\mc A}:\qcoh(G,X) \ra \qcoh(\mc A,X)$ is an equivalence of categories.

\end{prop}

\begin{proof}
We can first compose with the equivalence
$\pi^* :\qcoh Y \ra \qcoh(G,X)$, and by Proposition
\ref{EnrRightAdj} suppose that $\mc A$ has only one object $I$.
Define $\mc E = \pi_*^G (s_X^* I)$. Then the composite
functor is given by : $\mc G \ra \mc G \otimes_{\mc O_Y} \mc
E^{\vee}$. But we have also a functor $\mc F \ra \mc F \otimes_{\mc
  A_X} \mc E$ in the opposite direction, and natural units and
counits. To check that these are isomorphisms is a local problem, and
so is deduced by classical Morita theory, because by descent theory
$\mc E$ is locally free, hence is locally a progenerator.

\end{proof}

\begin{cor}
\label{FreeKtheory}
If the action of $G$ on $X$ is free then for each $i$ we have a natural
isomorphism $K_i(Y) \simeq K_i(\mc A,X)$.
\end{cor}

\begin{proof}
One checks that the equivalence of categories of Proposition
\ref{free=>equivalence} preserves coherence.
\end{proof}

\subsection{Euler characteristics of $\mc A$-sheaves}

Let $X$ be a proper $k$-scheme, endowed with the action of a finite
group $G$, and $\mc A$ be a full subcategory of $ k[G]\mod$,
admitting a finite set of additive generators.
Denoting by $s_X : X \ra {\rm spec}\, k$ the structure
morphism, we get a modular Euler characteristic $\chi(\mc A,\cdot) =
({s_X})_{\mc A} : K_0(\mc A, X) \ra  K_0(\mod \mc A)$ defined by
the usual formula $\chi(\mc A, \mc F)= \sum_{i\geq 0} (-1)^i [H^i(X,\mc F)]$.

For a sheaf $\mc G$ on $Y$, we denote as usual by $\chi(\mc G)$ its
ordinary Euler characteristic.

\begin{lem}
\label{red}
Let $\mc A$ be a projectively complete
full subcategory of $ k[G]\mod$,
admitting a finite set of additive generators.
Fix $\mc S$ a skeleton of the subcategory $\mc I$
of indecomposables of $\mc A$, and let $S$ be
the underlying finite set.
Let also $\mc F$ be an $\mc A$-sheaf on $X$.
Then have in $K_0(\mod \mc A)$ :
$$\chi(\mc A, \mc F)=\sum_{I\in S}\chi(\mc F(I))[S_I]$$
\end{lem}

\begin{proof}
Follows directly from the fact that $({s_X})_{\mc A}$ commutes with
restriction along $\langle I \rangle\ra \mc A$ and Lemma \ref{expl}.
\end{proof}

\begin{lem}
\label{red2} Suppose moreover that $\mc A$ is stable by radical
and duality, that duality and radical commute, and that $\mc F$
preserves monomorphisms. Then

$$\chi(\mc A, \mc F)=\sum_{I\in S}\sum_{i\geq 0}\chi({\rm gr}_i\mc F(I))[S_I]$$

\end{lem}

\begin{proof}
Since $\mc A$ admits a finite set of additive generators, the
canonical filtration of $\mc F$ (see \S \ref{CanFiltr}) is finite. Hence
$[\mc F]=\sum_{i\geq 0}[{\rm gr}_i \mc F]$ in $K_0(\mc A,X)$, and the
formula follows.
\end{proof}


\section{Applications}

\subsection{Symmetry principle}

We can prove in our context a formula given by Ellingsrud and L\o nsted (see
\cite{EL}, Theorem 2.4), following the proof of these authors :

\begin{prop}
\label{sym}
Let $X$ be a proper $k$-scheme, endowed with a \emph{free}
action of a finite group $G$,
and $\mc A$ be a full subcategory of $ k[G]\mod$,
admitting a finite set of additive generators, and containing
$k[G]$, the free object of rank $1$.
Then for each $G$-sheaf $\mc F$ on $X$ we have the equality in
$K_0(\mod\mc A)$ :

$$\chi(\mc A, \udl {\mc F})=\chi(\pi_*^G \mc F) \left[ \udl{k[G]}
\right] $$

\end{prop}

\begin{proof}
We can of course suppose that $\mc A$ is projectively complete.
  First, using Lemma \ref{expl}, it is easily seen that $\left[ \udl{k[G]}
\right] = \sum_{I\in\mc S}\dim _k I[S_I]$. Moreover Lemma \ref{red}
says that $\chi(\mc A,\udl {\mc F})=\sum_{I\in\mc S}\chi(\pi_*^G (I^{\vee}
\otimes_k \mc F ))[S_I]$. But since the action is free, $\pi $ is
\'etale, and for any sheaf $\mc G$ on $Y$, $\chi(\pi^* \mc G) =\# G
\chi(\mc G)$ (see \cite{Mum}), and the Proposition follows.
\end{proof}

Note that we recover the
formula of Ellingsrud and L\o nsted by evaluating at $k[G]$, in other
words, we have lifted the classical formula along
$K_0(\mod \mc A) \ra K_0(k[G] \mod)$.

\subsubsection{An example of computation}

Let $X$ be a smooth projective curve over an algebraically closed
field $k$ (i.e. a $1$-dimensional integral scheme, which is proper
over ${\rm spec}\, k$ and regular), endowed with a free action of a
finite group $G$. Let moreover $\mc A$ be a full subcategory of $ k[G]\mod$,
admitting a finite set of additive generators, and containing
$k[G]$, the free object of rank $1$.

We give some more details on the
cohomology of the sheaf of differentials on $X$.

For each object $V$ of $\mc A$, we choose a projective hull
$P(V)$, and then define $\Omega(V) = \ker P(V) \thra V$, and
inductively for $i\geq 0$ : $\Omega^{i+1}(V)= \Omega(\Omega^i(V))$,
$\Omega^0(V)=V$ (the notation $\Omega$ is unfortunate in our context,
but since it seems to be used in modular representation theory, we
keep it, hoping not to confuse the reader).

\begin{prop}
\label{diff}
$$H^0(X, \Omega_X) \oplus P(\Omega k) \simeq \Omega^2 k \oplus
k[G]^{\oplus g_Y -1} \oplus P(k)$$
\end{prop}

\begin{proof}
This is a special case of \cite{Kani}, Theorem 2.
\end{proof}

\begin{prop}
There is equality in $K_0(\mod\mc A)$ :
$$\left[ H^1(X, \udl {\Omega_X} ) \right]=\left[ \udl{\Omega^2 k}
\right] -\left[ \udl{P(\Omega k)}  \right] +\left[ \udl{P(k)} \right]$$
\end{prop}

\begin{proof}

Proposition \ref{diff} says that $ H^0(X, \Omega_X )$
is projectively equivalent to $\Omega^2 k$, so the comparison
Proposition \ref{comp} gives a five-terms exact sequence :

$$ 0 \ra {\rm Ext}^1_{k[G]}(\cdot, \Omega^2 k)
\ra H^1(X, \udl {\Omega_X} ) \ra \udl{k} \ra {\rm Ext}^2_{k[G]}(\cdot,
\Omega^2 k)\ra 0$$
Now, the long exact sequences associated to the short exact sequences :

$$0 \ra \Omega k \ra P(k) \ra k \ra 0$$
and

$$0 \ra \Omega^2 k \ra P(\Omega k) \ra \Omega k \ra 0$$
allow to conclude quickly.

\end{proof}

In our opinion, the above Proposition sheds some light on
the appearance of the $\Omega^2 k$ term in the Nakajima's
Proposition 3 of \cite{Naka},
and related formulas (like the one of Proposition \ref{diff}).
More precisely, Proposition \ref{diff} implies of course
that :

$$\left[ H^0(X, \udl {\Omega_X} ) \right]=\left[ \udl{\Omega^2 k}
\right] -\left[ \udl{P(\Omega k)}  \right] +\left[ \udl{P(k)} \right]
+\left[\udl{k[G]}^{\oplus g_Y -1}\right]
$$
and so $\udl {\Omega_X}$ verifies the symmetry principle in
$K_0(\mod \mc A)$, whereas $\Omega_X$ does not, in the sense that
$\left[ \udl{H^0(X,  \Omega_X )} \right]-\left[ \udl{H^1(X,
\Omega_X )} \right]$ is in general \emph{not} a multiple of
$\left[\udl{k[G]}
  \right]$ in $K_0(\mod \mc A)$.

\subsection{Galois modules on projective curves in positive
characteristic}

\subsubsection{Hypothesis}

\label{Hypo}

Following our previous paper \cite{Bor},
we show how to describe the
group $K_0(\mc A, X)$ when $X$ is a projective curve over an
algebraically closed field $k$. By projective curve we mean here a
$1$-dimensional integral scheme, which is proper over ${\rm spec}\, k$
and regular. We suppose that $X/k$ is endowed with a faithful action
of a finite group $G$. We will also suppose that $G$ acts with normal
stabilizers. We denote as usual by $\pi : X \ra Y=X/G$ the quotient.
$\mc A$ is a fixed subcategory of $k[G]\mod$ admitting a finite set of
additive generators.

\subsubsection{Additive structure of $K_0(\mc A, X)$}

For each $G$-invariant subset $U$, consider its complement $X'=X-U$,
endowed with the reduced structure, and denote by $f : X' \ra X$ the
corresponding closed immersion. Locally $X'$ is of the form
$G\times ^{G_P} P$, for a point $P$ of stabilizer $G_P$. Since $X'$ is
reduced, $P\simeq {\rm spec}\, k$, hence $G_P$ acts trivially on $P$.
Since moreover we make the hypothesis that $G_P$ is normal, we can
apply Proposition \ref{CritSurj} to deduce that $\mc A_X \ra
\tilde{f}_*\mc A_{X'}$ is an epimorphism. Now from Theorem \ref{Loc}
we get an exact sequence

$$\cdots \ra K_1(\mc A,U) \ra K_0(\mc A,X') \ra K_0(\mc A,X) \ra
K_0(\mc A,U) \ra 0$$

Since the category formed by all the $X'$, when $U$ varies between the
non empty $G$-invariant open subsets of $X$, is pseudofiltered, the
sequence remains exact after taking inductive limits on $X'$ (see
\cite{Schub} Theorem 14.6.6). Moreover the generic point $\xi$ of $X$
can be written as a $G$-scheme as
$$\xi = \lim_{\stackrel{\longleftarrow}{U}} U$$
We get as in \cite{Quil}, \S7, Proposition 2.2 that for each
nonnegative integer $i$

$$K_i(\mc A,\xi) \simeq \lim_{\stackrel{\longrightarrow}{U}} K_i(\mc A,U)$$

Moreover, since the action of $G$ on $X$ is faithful, $\xi$ is endowed
with a free action of $G$, with quotient $\eta$, the generic point of $Y$.
So according to Corollary \ref{FreeKtheory}, we have for each
nonnegative integer $i$ that $K_i(\eta)\simeq K_i(\mc A,\xi)$. As well
known, $K_0(\eta) \simeq \mathbb{Z}$, and $K_1(\eta)\simeq R(Y)^*$, where
$R(Y)$ is the function field of $Y$.

\begin{defi}
  (i) The {\emph group of $\mc A$-cycles on $X$}, denoted by $Z_0(\mc A, X)$, is
  by definition
  $$Z_0(\mc A, X)=\lim_{\stackrel{\longrightarrow}{X'}} K_0(\mc
  A,X')$$
  where the limit is taken on all the reduced strict closed $G$-subschemes
  of $X$.

  (ii) The {\emph group of classes of $\mc A$-cycles on $X$
    for the rational equivalence}, denoted by $A_0(\mc A, X)$,
  is by definition the cokernel of the
    canonical morphism $R(Y)^* \ra Z_0(\mc A, X)$ defined by the
    connection morphisms in the long exact sequences of $K$-theory.

   (iii) We denote by $\gamma :  A_0(\mc A, X) \ra K_0(\mc A,X)$ and
   ${\rm rk} :  K_0(\mc A,X) \ra \mathbb{Z}$ the canonical morphisms.
\end{defi}

\begin{thm}
\label{Mad}
Let $X$ be a projective curve over an algebraically closed field $k$,
endowed with the faithful action of a finite group $G$, acting with
normal stabilizers. Let moreover $\mc A$ be a subcategory of
$k[G]\mod$ admitting a finite set of additive generators, and
containing the free object $k[G]$.
Then the following morphism :
\begin{equation*}
 \xymatrix @R=0pt {
  \phi : \mathbb{Z} \oplus A_0(\mc A,X) \ar[r] & K_0(\mc A,X) \\
 (r,D) \ar[r] & r[\udl{\mc O_X}]+ \gamma(D) }
\end{equation*}
is an isomorphism.
 \end{thm}

\begin{proof}
As already seen, we have an exact sequence
$$ 0 \ra A_0(\mc A,X) \ra  K_0(\mc A,X) \ra \mathbb{Z} \ra 0$$
where the last map is the rank map. But this one is given by $ \mc F \ra
{\rm rk} (\mc F(k[G]))$, hence $r \ra r[\udl{\mc O_X}]$ is clearly a section.
\end{proof}

\begin{defi}
We denote by $c_1 : K_0(\mc A,X) \ra A_0(\mc A,X)$, and call
\emph {first Chern class}, the morphism composed of the inverse
$\phi^{-1} : K_0(\mc A,X) \ra \mathbb{Z} \oplus A_0(\mc A,X) $ of the
isomorphism of Theorem \ref{Mad}, followed by the second projection
$ \mathbb{Z} \oplus A_0(\mc A,X)\ra A_0(\mc A,X)$.
\end{defi}

\begin{lem}
The morphism ${\rm deg}_{\mc A}: Z_0(\mc A,X) \ra K_0(\mod \mc A)$
corresponding to the cone $((s_{X'})_{\mc A}: K_0(\mc A,X') \ra
K_0(\mod \mc A))_{X'}$ is trivial on the image of $R(Y)^* \ra  Z_0(\mc
A,X)$. We denote also by ${\rm deg}_{\mc A} : A_0(\mc A,X) \ra K_0(\mod \mc
A)$ the induced morphism.

\end{lem}

\begin{proof}
Let $\langle k \rangle$ denote the full subcategory of $k[G]\mod$ containing only
the trivial representation. Fix an object $V$ of $\mc A$,
corresponding to a unique $k$-linear functor $\langle k \rangle \ra \mc A$.
According to \cite{Quil}, \S5, Theorem 5,
the localization sequence is functorial for the restriction along this
functor, hence we obtain a commutative diagram :

\xymatrix {
 &&& K_1(\mc A, \xi) \ar[r] \ar[d] &  Z_0(\mc A,X) \ar[r]^{{\rm deg}_{\mc
  A}} \ar[d]  & K_0(\mc A,{\rm spec} \; k)\ar[d] \\
&&& K_1(\langle k \rangle, \xi) \ar[r]  &  Z_0(\mc \langle k \rangle,X) \ar[r]^{{\rm deg}_{\mc
 \langle k \rangle }}  & K_0(\langle k \rangle,{\rm spec} \; k)
}

But the bottom line is identified with $K_1(\eta) \ra Z_0(Y) \ra \mathbb{Z}$,
and \cite{Quil}, \S 7, Lemma 5.16 implies that the first map sends a
function to its divisor. Since $Y$ is also a projective curve,
the bottom line is thus a complex, and we can conclude from Lemma
\ref{expl} (note that by Morita invariance, we can reduce to the case
where $\mc A$ is projectively complete).
\end{proof}

\begin{cor}

\label{notsoexplicit}
Suppose the hypothesis of Theorem \ref{Mad} are verified.

Then for any coherent $\mc A$-sheaf $\mc F$ on $X$ we have in $K_0(\mod \mc
A)$ :

$$\chi(\mc A, \mc F) = {\rm rk} \mc F\; \chi(\mc A, \udl{\mc O_X})+
{\rm deg}_{\mc A} c_1(\mc F)$$

\end{cor}

\begin{proof}
According to Theorem \ref{Mad} we have
$[\mc F]= {\rm rk} \mc F\;[\udl{\mc O_X}]+c_1(\mc F)$ in $K_0( \mc
A,X)$, and the formula is obtained by pushing along $s_X : X \ra {\rm
  spec}\, k$.
\end{proof}

\subsubsection{Global sections of invertible sheaves pulled back from
  the quotient}

\begin{lem}
\label{vanishing}
There exists an integer $N$ such that
$$\forall \mc M \in {\rm Pic}\, Y \;\;\;
  \deg \mc M \geq N \Longrightarrow H^1(X, \udl{\mc \pi^*\mc M})=0$$
\end{lem}

 \begin{proof}
For a $k[G]$-module $V$, we have that
$H^1(X, \udl{\mc \pi^*\mc M})(V)=H^1(Y,\pi_*^G( \pi^*\mc M \otimes_{\mc
  O_X} s_X(V^{\vee})))= H^1(Y,\mc M \otimes_{\mc
  O_Y}\pi_*^G(  s_X(V^{\vee}))) $. The sheaf  $s_X(V^{\vee})$ is
locally free on the smooth curve $Y$, hence admits a filtration  whose
quotients are invertible, and so by Serre duality we can find $N_
  V$ so that for $ \deg \mc M \geq N_V$  we have $H^1(Y,\mc M \otimes_{\mc
  O_Y}\pi_*^G(  s_X(V^{\vee}))) =0$. Now we can choose $N= \sup_{V\in \mc
  S} N_V$, where $\mc S$ is a finite set of generators of $\mc A$.

 \end{proof} 

\begin{rem}
This result can be strengthened in the case of a cyclic group action,
see Proposition \ref{vanishing2}.
\end{rem}

Let $\mc G$ be an coherent sheaf on $Y$. To describe $c_1(\udl
{\pi^* \mc G})$, we first describe how to pull back cycles from the
quotient. 

\begin{lem}
For any $k[G]$-module $V$, the sheaf $\pi^G_* s_X^*(V)$ is locally
free of rank $\dim_k V$ on $Y$. 
\end{lem}

\begin{proof}
The sheaf $\pi^G_* s_X^*(V)$ is a subsheaf of $\pi_* s_X^*(V)$, hence
without torsion. But $Y$ is a smooth curve, hence the local rings at
its closed points are discrete valuation rings, and $\pi^G_* s_X^*(V)$
is in fact locally free. Moreover we have a monomorphism
$ \pi^* \pi^G_* s_X^*(V) \ra  s_X^*(V)$ with cokernel in the
ramification locus, which becomes an isomorphism at the generic point,
and so the two sheaves have same rank.  
\end{proof}

This Lemma allows to pull-back sheaves from
the quotient at the level of $K$-theory~: 

\begin{defi}
\

  (i) We denote by $\pi^{\mc A} : \coh \, Y \ra \coh(\mc A,X)$ the exact
  functor given on objects by $\mc G \ra \udl{\mc O_X} \otimes_{\mc
  O_Y} \mc G$, and also by  $\pi^{\mc A} : K_0(Y) \ra K_0(\mc A,X)$ the
  induced map.

  (ii) We define $\pi^{\mc A} : A_0(Y) \ra A_0(\mc A,X)$ as the unique
  morphism making the following diagram commute :

  \xymatrix{
   &&&  A_0(\mc A,X) \ar[r]^{\gamma} &  K_0(\mc A,X) \\
   &&& A_0(Y)\ar[u]^{\pi^{\mc A}} \ar[r]^{\gamma} & K_0(Y)\ar[u]^{\pi^{\mc A}}
    }
\end{defi}

\begin{lem}
\label{pullbackdiv}
\

(i) There is a functorial isomorphism $\pi^{\mc A} \mc G \simeq
\udl{\pi^* \mc  G}$

(ii) For any coherent sheaf $\mc G$ on $Y$ we have that $c_1( \pi^{\mc
  A} \mc G) = \pi^{\mc A} (c_1 \mc G)$.
\end{lem}

\begin{proof}
(i) is no more than the projection formula. 

(ii) results from the definitions of $\pi^{\mc A}$ and $c_1$.
\end{proof}







\begin{lem}
\label{degpullback}
\

(i) $\forall \delta \in A_0(Y) \;\;\; \deg_{\mc A} \pi^{\mc A} \delta =
  \deg \delta \cdot [\udl {k[G]}]$

(ii) $\forall \mc G \in \coh Y \;\;\;\deg_{\mc A} c_1(\udl {\pi^* \mc G}) =
\deg (c_1 \mc G) \cdot [\udl {k[G]}]$
\end{lem}

\begin{proof}
(ii) is a consequence from (i) and Lemma \ref{pullbackdiv}.
 For (i), we can choose $\delta=Q$, a reduced closed point.
 Let $X' = (\pi^{-1} Q)^{red}$ be its  reduced preimage by $\pi$,
 $i : X' \ra X$ denote the closed immersion, and
 $j : Q \ra Y$ the corresponding morphism between the quotients.
 By d\'evissage, we can choose $\pi^{\mc A}\delta$ to be represented
 by an element of $K_0(\mc A,X')$, which we denote the same way.
 By definition of $c_1$, we have that
 $i_{\mc A}(\pi^{\mc A}\delta)=[\pi^{\mc A}j_*\mc O_Q]$ in $K_0(\mc
 A,X)$. It follows that $\deg_{\mc A}(\pi^{\mc A}\delta)=(s_{X'})_{\mc
 A}((\pi^{\mc A}\delta))= (s_{X})_{\mc A} [\pi^{\mc A}j_*\mc O_Q]
 = [(s_{X})_{\mc A} \pi^{\mc A}j_*\mc O_Q]$. Note that $[\pi^{\mc A}  j_* \mc
   O _Q]= [\udl {\mc O_X} \otimes_{O_Y} j_* \mc O _Q]=   [V\ra \pi^G_*
 s_X^*(V^{\vee})\otimes _{O_Y}j_* \mc O _Q]=[V \ra   j_*j^*\pi^G_*
 s_X^*(V^{\vee})]$. Finally we get that $\deg_{\mc A}(\pi^{\mc
 A}\delta)=[V\ra H^0(Q,j^*\pi^G_*
 s_X^*(V^{\vee})]$. Since $\dim_k H^0(Q,j^*\pi^G_*
 s_X^*(V^{\vee})= \dim_k V$, we can conclude thanks to Lemma 
\ref{expl} (again, by Morita invariance, we can reduce to the case
where $\mc A$ is projectively complete).

\end{proof} 

\begin{rem}
Instead of  $ \deg_{\mc A} \pi^{\mc A} \delta$, it would be better to
describe $\pi^{\mc A} \delta$ itself. It is quite clear that it is a
cycle with ``free coefficients in the inertia'', as it can be probably
be deduced from Lemma \ref{degpullback}.
\end{rem}

\begin{thm}
\label{pullbackinv}
Let $X$ be a projective curve over an algebraically closed field $k$
of positive characteristic $p$, endowed with the faithful action of
finite group $G$. Suppose that $G$ has cyclic $p$-Sylows  and
acts with normal stabilizers.
Denote by $\pi : X \ra Y=X/G$ the quotient map.
There exists an integer $N$ such that for all invertible sheaves $\mc
M$, $\mc N$ on $Y$ such that $\deg \mc M \geq \deg  \mc N \geq N$ the
following isomorphism of $k[G]$-modules holds :

$$ H^0(X, \pi^* \mc M) \simeq H^0(X, \pi^*\mc N) \oplus k[G]^{\oplus
  \deg \mc M - \deg \mc N}$$

\end{thm}

\begin{proof}
We choose $\mc A=k[G]\mod$.
According to Corollary \ref{notsoexplicit}, Lemma \ref{vanishing} and
Lemma \ref{degpullback} (ii), we have for
$\deg \mc M \geq N$ that $[\udl{ H^0(X,\pi^*\mc M)}]=
\chi(\mc A, \udl{\mc O_X})+
\deg \mc M [\udl {k[G]}]$ in $K_0(\mod \mc A)$, and similarly for $\mc
N$, and we can conclude thanks to Theorem \ref{Main} (iii).

\end{proof}

\subsubsection{Explicit expression for the action of a cyclic group}

\label{Cyclic}

To the hypothesis of \S \ref{Hypo}, we add the fact that $G$ is cyclic
of order $n$, and we fix a generator $\sigma$. Moreover we choose as
ring with several objects $\mc A = k[G] \, \mod$, the whole category
of all $k[G]$-modules of finite type. We follow the notations given in
\S \ref{CanFiltr}.

So by definition, for each
coherent kernel-preserving
$\mc A$-sheaf $\mc F$ on $X$, we denote by  ${\rm gr}_0 \mc
F $ the $\mc A$-sheaf on $X$ defined on indecomposables by the
existence of $n$ exact sequences of sheaves on $Y$ :

$$ 0 \ra \mc F(\psi \otimes V_{j-1}) \ra \mc F(\psi \otimes V_j) \ra
{\rm gr}_0 \mc F (\psi \otimes V_j)\ra 0$$

The reader who is just concerned by the cyclic case may
skip \S \ref{CanFiltr} and consider this as a definition of ${\rm gr}_0
\mc F$, since this is all what we need here.

\begin{lem}
\label{explicit}

For any coherent $G$-sheaf $\mc G$ on $X$ :
$$ \chi (\mc A, \udl {\mc G})
=\sum_{\psi \in \widehat G} \sum_{j=1}^{p^v} \sum_{i=1}^{j}
\chi(\emph{gr}_0 \udl {\mc G}(\psi \otimes V_{i}))[S_{\psi \otimes V_{j}}]$$

\end{lem}

\begin{proof}
This follows from Lemma \ref{red2}, by noticing that in the
particular case of a cyclic group $G$ we have ${\rm gr}_i {\mc
  F} \simeq R^i{\rm gr}_0  {\mc F}$. Alternatively, this is direct
from the above exact sequences and Lemma \ref{red}.
\end{proof}

From now on, we will consider an invertible $G$-sheaf $\mc L$ on
$X$, and show how to determine explicitly the structure of
$k[G]$-module of its global sections $H^0(X, \mc L)$ when $\deg
\mc L > 2g_X-2$.

To do this, we will of course use Lemma \ref{explicit}, and show
that we can describe explicitly the invertible sheaves ${\rm gr}_0
\udl {\mc L} (\psi \otimes V_j)$ in terms of a $G$-invariant
divisor $D$ such that $\mc L \simeq \mc L_X (D)$, of the
ramification data, of the genus of $Y$, and of some rational
functions on $X$.

Notice that the formula
$$\pi_*^G \mc L_X (D) \simeq \mc L_Y ([\pi_* D/\#G ])$$
(where $[\cdots ]$ is the integral part of the divisor, taken
coefficient by coefficient) provides an explicit
description of $\pi_*^G \mc L$ (see \cite{Kani}, proof of Proposition 3).

First we reduce to the case of a cyclic $p$-group. Remember that
$n=p^v a$, with $a$ prime to $p$, and let $H$ be the subgroup of
$G$ of order $a$, and $P=G/H$, so that we have a tower :

\begin{equation}
  \label{tower}
\xymatrix{
     X \ar[dd]_G^{\pi} \ar[dr]^H_{\alpha} &  \\
     &  Z \ar[dl]^P_{\beta} \\
     Y  &   }
\end{equation}

\begin{lem}
For each $\psi$ in the character group $\widehat G$, there exists a
nonzero function $f_\psi$ in the function field $R(X)$ of $X$ such
that for each invertible $G$-sheaf $\mc L$ on $X$ we have :

$${\rm gr}_0^G \udl {\mc L} (\psi \otimes V_j) \simeq
  {\rm gr}_0^P \udl {\alpha_*^H  {\mc L}(f_\psi)} (V_j) $$
\end{lem}

\begin{proof}
From Kummer theory, we get for each $\psi$ in $\widehat G$ a non
zero function $f_\psi$ such that $\mc L_X((f_\psi)) \simeq
s_X^*(\psi)$. Thus the result follows from the exactness of $\alpha_*^H$.

\end{proof}

Next we reduce to the case of a cyclic group of order $p$. For this,
suppose that we start from a group $G$ of order $p^v$, with $v \geq
2$, and let $H$ be the subgroup of order $p^{v-1}$, and $P=G/H$. We
keep the notations of diagram \ref{tower} in this context.

\begin{prop}

\label{recur}
Let $1\leq j \leq p^{v}$ an integer, and write :
$ j=(l-1)p+j'$
with $1\leq l \leq p^{v-1}$ and $1\leq j' \leq p$.

Then we have an isomorphism of invertible sheaves on $Y$ :

$$ {\rm gr}_0^G  \udl {\mc L}(V_j) \simeq
{\rm gr}_0^P \udl{{\rm gr}_0^H\udl{\mc L}(V_l)}(V_{j'})$$
\end{prop}

\begin{proof}

The proof requires the following Lemmas.

\begin{lem}

\label{decomp}
$${V_j}_{|H} \simeq V_l^{\oplus j'} \oplus V_{l-1}^{\oplus p-j'}$$
\end{lem}

\begin{proof}
$1$, $\sigma-1$, ... , $(\sigma-1)^{j'-1}$ generate Jordan blocks of
  size $l$, while $(\sigma-1)^{j'}$, $(\sigma-1)^{j'+1}$, ... ,
  $(\sigma-1)^{p-1}$, generate Jordan blocks of size $l-1$.
\end{proof}

\begin{lem}

\label{onZ}
There is an exact sequence of $P$-sheaves on $Z$ :

$$ 0 \ra \alpha_*^H(\mc L \otimes V_{(l-1)p}) \ra
\alpha_*^H(\mc L \otimes V_j) \ra V_{j'} \otimes  {\rm gr}_0^H(\mc
L)(V_l) \ra 0$$
\end{lem}

\begin{proof}
First ignoring the action of $P$, we can use the Lemma \ref{decomp} to
show the existence of an exact sequence of sheaves on $Z$ :

$ 0 \ra \alpha_*^H(\mc L \otimes V_{(l-1)p}) \ra
\alpha_*^H(\mc L \otimes V_j) \ra ({\rm gr}_0^H(\mc
L)(V_l))^{\oplus j'} \ra 0$

But then from the exact sequence of $k[G]$-modules
$$0 \ra V_{(l-1)p} \ra V_j \ra V_{j'} \ra 0$$ one sees that $\sigma -1$
acts transitively
on the direct summands of $({\rm gr}_0^H(\mc L)(V_l))^{\oplus j'}$, hence the
result.
\end{proof}

We will suppose that $j'\geq 2$, the case $j'=1$ being analog, and easier.
We have the following commutative diagram of $P$-sheaves on $Z$ :

\xymatrix{
&& 0 & 0 & \\
&& {\rm gr}_0^H(\mc L)(V_l)\ar[u]\ar@{=}[r] & {\rm gr}_0^H(\mc L)(V_l)\ar[u]& \\
0  \ar[r] & \alpha_*^H(\mc L \otimes V_{(l-1)p}) \ar[r]
&   \alpha_*^H(\mc L \otimes V_j) \ar[r]\ar[u] & V_{j'} \otimes  {\rm gr}_0^H(\mc
L)(V_l)\ar[r]\ar[u] & 0 \\
0 \ar[r] & \alpha_*^H(\mc L \otimes V_{(l-1)p})\ar[r] \ar@{=}[u]&\alpha_*^H(\mc L \otimes
V_{j-1}) \ar[r]\ar[u] & V_{j'-1} \otimes  {\rm gr}_0^H(\mc
L)(V_l)\ar[r]\ar[u]& 0 \\
&&0 \ar[u] &0 \ar[u]&}

So the Proposition follows now from a diagram chase thanks to the
following fact~:

\begin{lem}
$$R^1\beta_*^P(\alpha_*^H(\mc L \otimes V_{(l-1)p}))=0$$
\end{lem}

\begin{proof}
One sees as in Lemma \ref{onZ} that
$\alpha_*^H(\mc L \otimes V_{(l-1)p}) \simeq V_p \otimes
\alpha_*^H(\mc L \otimes V_{l-1})$, and since $V_p=k[P]$, this is enough.
\end{proof}

\end{proof}

So we are reduced to the case of a cyclic $p$-group, which was solved
by S.Nakajima. We give a translation of \cite{Naka}, Theorem 1 in our
context, and since our version is slightly stronger, we give a sketch
of a proof.

\begin{thm}[Nakajima]
\label{incre}
Let $X$ be a projective curve endowed with a faithful action of
$G=\mathbb Z/p$, and $D$ a $G$-invariant divisor on $X$. Write
$D=\pi^* \delta +\sum_{P\in X_{ram}}n_P \cdot P$, where $\delta$ is a
divisor on $Y$ so that ${\rm supp} \;\pi^*\delta \cap X_{ram} =
\varnothing$. For each $P$ in $X_{ram}$, let moreover $N_P$ be the
integer defined by $N_P+1=v_P(\sigma u_P -u_P)$, where $\sigma$ is a
generator of $G$, and $u_P$ an uniformizer at $P$. Then for each
integer $1 \leq j \leq p$ :

$$ {\rm gr}_0^G  \udl {\mc L_X(D)}(V_j) \simeq \mc L_Y (\delta +
\sum_{P\in X_{ram}} [\frac{n_P-(j-1)N_P}{p}]\cdot \pi_* P)$$
\end{thm}

\begin{proof}

A local analysis shows (see Lemma \ref{local}) that the monomorphism
$\pi^*  {\rm gr}_0^G  \udl {\mc L_X(D)}(V_j) \ra \mc L_X(D)$
factorizes trough
$\pi^*  {\rm gr}_0^G  \udl {\mc L_X(D)}(V_j) \ra \mc L_X(D-
\sum_{P\in X_{ram}}(j-1)N_P\cdot P)$, hence by applying $\pi_*^G$ we
get a monomorphism $ {\rm gr}_0^G  \udl {\mc L_X(D)}(V_j) \ra
\mc L_Y (\delta + \sum_{P\in X_{ram}} [\frac{n_P-(j-1)N_P}{p}]\cdot
\pi_* P)$. To show that this is an isomorphism is a local problem at
$X_{ram}$, so by adding eventually to $D$ a divisor of the form $\pi^*
\gamma$ we may suppose that $\deg D > 2g_X-2$. But then \cite{Naka},
Lemma 4, shows that the two sheaves have the same space of global
sections, hence they must be isomorphic.
\end{proof}

\begin{lem}
\label{local}
With the notations of Theorem \ref{incre} we have :

$$\pi^*  {\rm gr}_0^G  \udl {\mc L_X(D)}(V_j) \subset \mc L_X(D-
\sum_{P\in X_{ram}}(j-1)N_P\cdot P)$$

\end{lem}

\begin{proof}
This is a local problem, so we can check the inclusion
on the completions of the local rings of the closed points of $X$. Let
$P$ be such a point, that we can suppose in $X_{ram}$, $Q=\pi P$,
$u_P$ an uniformizer at $P$, $v_Q$ a uniformizer at $Q$. Localizing at
$P$ the commutative diagram :

$$
\xymatrix{
  0 \ar[r] &\mc L_X(D) \otimes V_{j-1}  \ar[r] &\mc L_X(D) \otimes V_j
  \ar[r] & \mc L_X(D) \ar[r] & 0 \\
  0 \ar[r] & \pi^*\pi_*^G (\mc L_X(D) \otimes V_{j-1})  \ar[r] \ar[u] &
  \pi^*\pi_*^G(\mc L_X(D) \otimes V_j)
  \ar[r] \ar[u] & \pi^*  ({\rm gr}_0^G  \udl
  {\mc L_X(D)}(V_j)) \ar[r] \ar[u]& 0
  }
$$

we get the following commutative diagram :

$$
\xymatrix{
\cdots \ar[r]
& V_j\otimes_k u_P^{-n_P} k[[u_P]] \ar[r]^{\mu} & u_P^{-n_P} k[[u_P]] \ar[r] & 0\\
\cdots \ar[r]
& k[[u_P]]\otimes_{k[[v_Q]]}(V_j\otimes_k u_P^{-n_P} k[[u_P]])^G
\ar[r] \ar[u] & k[[u_P]]\otimes_{k[[v_Q]]}(\widehat{{\rm gr}_0^G  \udl
  {\mc L_X(D)}(V_j)_Q}) \ar[r]\ar[u]^{\nu} & 0}$$

What we have to show is that the image of $\nu$ is contained in the
ideal $u_P^{-n_P+(j-1)N_P} k[[u_P]]$. Choose $\alpha$ in
$V_j\otimes_k u_P^{-n_P} k[[u_P]]$. Since $V_j = k[G]/(\sigma
-1)^j$, we can write $\alpha=\sum_{i=0}^{j-1}(\sigma -1)^i a_i$, with
$a_i \in u_P^{-n_P} k[[u_P]]$. Then $\mu(\alpha)=a_0$, and if $\sigma
\alpha = \alpha$, then $a_0=(\sigma^{-1}-1)^{j-1}a_{j-1}$. But now we
can conclude, since the definition of $N_P$ implies that for any $x$
in $k[[u_P]]$, $v_P((\sigma^{-1}-1)x) \geq v_P(x) + N_P$.
\end{proof}

Note moreover that in a relative situation like those appearing in
various d\'evissage steps, the a priori nonequivariant isomorphisms are in
fact automatically equivariant~: indeed when $G$ is a $p$-group acting
on a projective $k$-scheme $X$, there is a most a structure of
$G$-sheaf on a given invertible sheaf. So Proposition \ref{recur}
allows to apply Theorem \ref{incre} recursively, to finally have
an explicit expression of ${\rm gr}_0^G  \udl {\mc
  L_X(D)}(V_j)$, i.e. to represent this invertible sheaf on $Y$ by a
divisor.

Once these sheaves are computed, we return to the case of an
arbitrary cyclic group, and show how we can make use of Lemma
\ref{explicit}.

\begin{prop}
\label{vanishing2}

  Suppose that $G$ is any cyclic group.
If $\deg \mc L > 2g_X-2$ then $H^1(X, \udl{\mc L})=0$.
\end{prop}

\begin{proof}
This is equivalent to the fact that for each object $V$ of $\mc A$, the
sheaf $\udl{\mc L}(V)$ on $Y$ is acyclic, and this can be checked on
indecomposables. But, as
one sees by induction, it is enough to show that for each $\psi$ in
$\widehat G$, and $1\leq j \leq p^v$, the invertible sheaf
${\rm gr}_0^G  \udl {\mc L}(\psi \otimes V_j)$ on $Y$ is acyclic.
To show this, we follow the same reduction steps. The reduction to the
case of a cyclic $p$-group is an immediate computation using the Hurwitz
formula in its basic (tame) form. The Proposition \ref{recur} allows
then to reduce to the case of a cyclic group of order $p$. But then
Theorem \ref{incre} and the computation in the proof of \cite{Naka},
Lemma 3, allows to conclude.
\end{proof}

Now Lemma \ref{explicit} provides, for an invertible sheaf $\mc L$
such that $\deg \mc L > 2g_X-2$, an expression in $K_0(\mod \mc A)$:

\begin{equation}
\label{Exp}
  [\udl{H^0(X,\mc L)}]
=\sum_{\psi \in \widehat G} \sum_{j=1}^{p^v} \sum_{i=1}^{j}
\chi({\rm gr}_0 \udl {\mc L}(\psi \otimes V_{i}))[S_{\psi \otimes
    V_{j}}]
\end{equation}
where the integers $\chi({\rm gr}_0 \udl {\mc G}(\psi \otimes
V_{j-i}))$ are given by the previous description and the classical
Riemann-Roch formula. Theorem \ref{Main}, (ii), implies that this
characterizes fully $H^0(X,\mc L)$ as a $k[G]$-module. But we can be
more explicit by applying Theorem \ref{Main}, (iii). For this we have to use
the basis
$\{ [\udl{\psi \otimes V_j}] \; / \; \psi \in \widehat G, \; 1 \leq j \leq
p^v \}$ of $K_0(\mod \mc A)$ rather that the basis
$\{ S_{\psi \otimes V_j} \; / \; \psi \in \widehat G, \; 1 \leq j \leq
p^v \}$.

Clearly, we can find the base change matrix block by block (i.e.
here character by character) and thus fix a $\psi$ in $\widehat
G$. Let $A$ be the Cartan matrix, i.e. the matrix whose $j$-th
column  is formed by the coefficients of $[\udl{\psi \otimes
V_j}]$ in the basis $\{ S_{\psi \otimes V_i}\; / \; 1 \leq i \leq
p^v \}$. Then Lemma \ref{expl} shows that the coefficients of $A$
are given by $\dim_k {\rm Hom}_{k[G]}(V_j,V_i)$, and this is
easily computed as being $\inf(i,j)$. This is independent of
$\psi$, so all the blocks of the base change matrix are equal, and
the inverse is easily computed as :

\begin{equation}
\label{Matrix}
  A^{-1}=\left(
     \raisebox{0.5\depth}{%
      \xymatrixcolsep{1ex}%
       \xymatrixrowsep{1ex}%
       \xymatrix{
         + 2 & -1 \ar @{-}[dddrrr] &
         0 \ar @{.}[dddrrr] \ar@{.}[rrr]
         & & & 0 \ar@{.}[ddd]\\
         -1 & +2 \ar @{-}[dddrrr] &&&& \\
         0 \ar@{.}[ddd] \ar @{.}[dddrrr] & -1 \ar @{-}[dddrrr] &&&&\\
         &&&& -1&0 \\
     &&&& +2&-1 \\
         0 \ar@{.}[rrr]& & &0 & -1 & +1
       }%
     }
   \right)
\end{equation}

\subsubsection{Explicit expression for the action of a cyclic $p$-group}
\label{formul}

\begin{defi}
Let $\pi : X\ra Y$ be a (generically) Galois cover of projective
curves over $k$ of group $G \simeq \mathbb Z/ p$, with generator
$\sigma$.

(i) For a ramified point $P$ of $X$, define $N_P$ as the integer
such that $N_P+1$ is the valuation at $P$ of $\sigma u_P-u_P$,
where $u_P$ is an uniformizer at $P$.

(ii) For $0\leq \alpha\leq p-1$ define a map $\pi_*^ \alpha:
Z_0(X) \ra Z_0(Y)$ between $0$-cycles groups by, for any divisor
$D$ on $X$, $$\pi_*^\alpha D = [ \frac{1}{p}\pi_*(D-\alpha\sum_{P
\in
  X_{ram}}N_P \cdot P)]$$ where $[\cdots]$ denotes the integral part of a
divisor, taken coefficient by coefficient.

\end{defi}

\begin{thm}
\label{formula}
 Suppose that $X$ is a projective curve over $k$ with a faithful
 action of $G \simeq \mathbb Z/ p^v  \mathbb Z$. For $1\leq n\leq v$
 let $X_n$ be the quotient curve of $X$ by the action of the subgroup
 of $G$ of order $p^n$, and $\pi_n : X_{n-1} \ra X_n$ be the canonical
 morphism. Let moreover $D$ be a $G$-invariant divisor on $X$, and
 $H^0(X,\mc L_X(D))\simeq \oplus_{j=1}^{p^v}V_j^{\oplus m_j}$ be the
 Krull-Schmidt decomposition of the global sections of $\mc L_X(D)$,
 where $V_j$ is the indecomposable $k[G]$-module of dimension $j$.
 Suppose $\deg D > 2g_X-2$. Then the integers $m_j$ are given by :
$$\left\{ \begin{array}{lr}

m_j=\deg({\pi_v}_*^{\alpha_0(j)} \cdots
{\pi_1}_*^{\alpha_{v-1}(j)}D)- \deg({\pi_v}_*^{\alpha_0(j+1)}
\cdots {\pi_1}_*^{\alpha_{v-1}(j+1)}D)
& {\rm  if}\;\; 1\leq j\leq p^v-1\\
m_{p^v}=1-g_{X_v} +  \deg({\pi_v}_*^{p-1} \cdots {\pi_1}_*^{p-1}D)
& \end{array}                      \right. $$

\noindent where  for $1\leq j\leq p^v$ the integers
${\alpha_0(j)},\cdots,{\alpha_{v-1}(j)}$ are the digits of the
$p$-adic writing of $j-1$ defined by
$j-1=\sum_{h=0}^{v-1}\alpha_h(j) p^h$ with $0\leq \alpha_h(j)\leq
p-1$.
\end{thm}

\begin{proof}

This is a direct consequence from the following intermediate
computation :

\begin{lem}

$$\left\{ \begin{array}{lr}
m_1=2a_1-a_2 & \\
m_j=-a_{j-1}+2a_j-a_{j+1} & {\rm  if}\;\; 2\leq j\leq p^v-1\\
m_{p^v}=a_{p^v}-a_{p^v-1} & \end{array}                      \right. $$

\noindent with, for  $1\leq j\leq p^v$,
$a_j=j(1-g_{X_v})+\sum_{i=1}^j \deg({\pi_v}_*^{\alpha_0(i)} \cdots
{\pi_1}_*^{\alpha_{v-1}(i)}D)$.
\end{lem}

\begin{proof}
  Set $\mc L=\mc L_X(D)$ and $X_v=Y$.
According to the equation \ref{Exp} we have
  $[\udl{H^0(X,\mc L)}]
= \sum_{j=1}^{p^v} \sum_{i=1}^{j}
\chi({\rm gr}_0 \udl {\mc L}( V_{i}))[S_{
    V_{j}}]$. Set $b_j=\sum_{i=1}^{j}
\chi({\rm gr}_0 \udl {\mc L}( V_{i}))$. If we show that for each $j$
  we have equality $b_j=a_j$, we are done, according to the base
  change matrix given in equation \ref{Matrix}. But the usual
  Riemann-Roch formula gives $b_j=j(1-g_Y)+\sum_{i=1}^{j} \deg ({\rm
    gr}_0 \udl {\mc L}( V_{i}))$. So the only thing which remains to
  be shown is ${\rm
    gr}_0 \udl {\mc L}( V_{i}) \simeq \mc L_Y({\pi_v}_*^{\alpha_0}
  \cdots {\pi_1}_*^{\alpha_{v-1}}D)$.
  This is easily done by induction : the first
  step is Theorem \ref{incre}, and the induction step is given by
  Proposition \ref{recur}.
\end{proof}

\end{proof}

\subsubsection{Noether criterion with parameter}

\begin{defi}
Let $G$ a finite group, $H$ a subgroup, $k$ an algebraically
closed field. A $k[G]$-module of finite type $V$ is said
\emph{relatively $H$-projective} if it is a direct summand of a
module induced from $H$.
\end{defi}

\begin{thm}
\label{Noether}
 Let $\pi : X\ra Y$ be a (generically) cyclic
Galois $p$-cover of projective curves over $k$ of group $G$, ${\rm
ram}\,\pi$ the largest ramification subgroup of $\pi$, and $H$ a
subgroup of $G$. Then the following assertions are equivalent :

(i) ${\rm ram}\,\pi \subset H$

(ii) $\forall \mc L \in {\rm Pic}_G X \;\;\; \deg \mc L > 2g_X-2
\Longrightarrow H^0(X, \mc L)$ is relatively $H$-projective.

(iii) $\exists \mc M \in {\rm Pic}\, Y$ so that $\# G \deg \mc M >
2g_X-2$ and $H^0(X, \pi^*\mc M)$ is relatively $H$-projective.

\end{thm}

\begin{proof}
We begin by some notations. Let $p^v$ (resp. $p^w$) be the order
of $G$ (resp. $H$). Fixing  a generator $\sigma$ of $G$, we denote
for $1\leq l\leq p^v$ the $k[G]$-indecomposable of dimension $l$
by $V_l^G=k[\sigma]/(\sigma -1)^l$, and similarly for $H$.

\begin{lem}
\label{Hproj}
 The relatively $H$-projective $k[G]$-modules are
exactly the direct sums of the indecomposables $V_{lp^{v-w}}^G$
for $1\leq l\leq p^w.$
\end{lem}

\begin{proof}
This results from the Krull-Schmidt Theorem and the fact that
${\rm Ind}_H^G V_l^H =V_{lp^{v-w}}^G$ for $1\leq l\leq p^w.$
\end{proof}

We will also use the notations and results of \S \ref{formul}.

We first show that (i) implies (ii). For this, note that according
to Theorem \ref{formula}, the multiplicity of $V_j^G$ in $H^0(X,
\mc L)$ is given for $1\leq j\leq p^v-1$ by :

$m_j=\deg({\pi_v}_*^{\alpha_0(j)} \cdots
{\pi_{w+1}}_*^{\alpha_{v-w-1}(j)}{\pi_w}_*^{\alpha_{v-w}(j)}
\cdots   {\pi_1}_*^{\alpha_{v-1}(j)}D)-\newline
\deg({\pi_v}_*^{\alpha_0(j+1)} \cdots
{\pi_{w+1}}_*^{\alpha_{v-w-1}(j+1)}{\pi_w}_*^{\alpha_{v-w}(j+1)}\cdots
{\pi_1}_*^{\alpha_{v-1}(j+1)}D) $

Because of the hypothesis that ${\rm ram}\,\pi \subset H$, the
coverings $\pi_u$ are \'etale for $u>w$, hence the morphisms
${\pi_u}_*^{\alpha}$ are independent of $\alpha$. To show (ii) we
can, according to the Lemma \ref{Hproj}, show that $m_j=0$ for
$p^{v-w} \nmid j$, which results from the previous remark and of :

\begin{lem}
Let $1 \leq j\leq p^v$ so that $p^{v-w} \nmid j$. Then $\forall u
\geq v-w \;\;\alpha_u(j)=\alpha_u(j+1)$.
\end{lem}

\begin{proof}
This is clear from the definition of $\alpha_u(j)$ as the $u+1$-th
digit in the $p$-adic writing of $j-1$.
\end{proof}

Since the fact that (ii) implies (iii) is trivial, all what
remains to be shown is that (iii) implies (i). For this suppose
that ${\rm ram}\,\pi \nsubseteq H$, that is, since the subgroups
of $G$ are totally ordered, $H \subsetneqq {\rm ram}\,\pi$. Let
$\mc M$ be an invertible sheaf on $Y$ such that $\# G \deg \mc M >
2g_X-2$. To show that (iii) is false, we can clearly suppose that
${\rm ram}\pi$ has order $p^{w+1}$. Choose $j=p^{v-w-1}$. An easy
computation shows that the coefficient $m_j$ of $V_j^G$ in $H^0(X,
\pi^*\mc M)$ is

$$-\sum_{P\in X_w^{ram}}[-N_P/p]$$

and since by hypothesis $X_w^{ram}\neq \emptyset$, this is non
zero. Hence $H^0(X, \pi^*\mc M)$ is non relatively $H$-projective,
as was to be shown.

\end{proof}

\bibliography{Arxiv3}
\end{document}